\documentclass[10pt,reqno]{amsart}

\usepackage{amsfonts,latexsym,amsmath,amssymb,amscd,amsmath}
\usepackage{xypic, epsfig}
\xyoption{all}

\theoremstyle{plain}
\newtheorem{theorem}{Theorem}[section]
\newtheorem{proposition}[theorem]{Proposition}
\newtheorem{corollary}[theorem]{Corollary}
\newtheorem{lemma}[theorem] {Lemma}

\theoremstyle{definition}
\newtheorem{definition}{Definition}[section]

\newtheorem{example}{Example}[section]

\theoremstyle{remark}
\newtheorem{remark}{Remark}[section]

\providecommand{\e}{\ar@{-}}
\providecommand{\de}{\ar@{.}}

\providecommand{\E}{\mathcal{E}}
\providecommand{\cI}{\mathcal{I}}
\providecommand{\calD}{\mathcal{D}}
\renewcommand{\Im}{\mathfrak{Im}}
\providecommand{\al}{\alpha}
\providecommand{\be}{\beta}
\providecommand{\la}{\lambda}
\providecommand{\bC}{\mathbb{C}}
\providecommand {\bR}{\mathbb R}

\providecommand {\DD} {\mathbf D}

\providecommand {\Ga} {\Gamma}
\providecommand {\ga} {\gamma}
\providecommand {\eps} {\epsilon}
\providecommand\pa {\partial}

\providecommand{\bN}{\mathbb N}
\providecommand{\HH}{\mathcal{H}}

\providecommand{\bZ}{\mathbb Z}

\providecommand {\A} {\mathcal A}

\providecommand{\NN}{\mathbb N}

\providecommand{\RR}{\mathbb{R}}
\providecommand{\ZZ}{\mathbb{Z}}
\providecommand{\CC}{\mathbb{C}}
\providecommand{\FS}{\mathcal{F}}

\renewcommand{\S}{Section }

\def\newop#1{\expandafter\def\csname #1\endcsname{\mathop{\rm
#1}\nolimits}}

\newop{slc}
\newop{sign}
\newop{mdeg}
\newop{supp}

\title[Multivariate P\'olya-Schur Classification Problems]
{Multivariate P\'olya-Schur Classification Problems in the Weyl Algebra}

\author[J.~Borcea]{Julius Borcea}
\address{J.~Borcea, Department of Mathematics, Stockholm University, 
SE-106 91 Stockholm, Sweden}
\author[P.~Br\"and\'en]{Petter Br\"and\'en}
\address{P.~Br\"and\'en, Department of Mathematics, Royal Institute of Technology, 
SE-100 44 Stockholm, Sweden, and,  Department of Mathematics, Stockholm University, 
SE-106 91 Stockholm, Sweden}
\email{pbranden@math.su.se}

\begin{document}

\begin{abstract} 
A multivariate polynomial is {\em stable} if it is nonvanishing whenever all 
variables have positive imaginary parts. We classify all linear partial 
differential operators in the Weyl algebra $\A_n$ that preserve stability. 
An important tool that we develop in the process is the higher dimensional 
generalization of P\'olya-Schur's notion of multiplier sequence. 
We characterize all multivariate multiplier sequences as 
well as those of finite order. 
Next, we establish a multivariate extension of the Cauchy-Poincar\'e 
interlacing theorem and prove a natural analog of 
the Lax conjecture for real stable polynomials in two variables. Using the 
latter we describe all operators in $\A_1$ that preserve univariate
hyperbolic polynomials by means of determinants and homogenized symbols. 
Our methods also yield homotopical properties for 
symbols of linear stability preservers and a duality theorem showing that an  
operator in $\A_n$ preserves stability if and only if its Fischer-Fock 
adjoint does. These are powerful multivariate extensions of the classical 
Hermite-Poulain-Jensen theorem, P\'olya's curve theorem and 
Schur-Mal\'o-Szeg\H{o} composition theorems. 
Examples and  applications to strict stability preservers 
are also discussed.
\end{abstract}

\maketitle

\tableofcontents

\section{Introduction and main results}\label{s-intro}

In their seminal 1914 paper \cite{PS} P\'olya and Schur characterized all
linear operators that are diagonal in the standard monomial basis of $\bC[z]$
and preserve the set of polynomials with only real zeros. Polynomials of this
type and linear transformations preserving them are of central interest in 
e.g.~entire function theory \cite{CC1,Le}: it is for instance well known that
the Riemann Hypothesis is equivalent to saying that 
$\xi\!\left(\frac{1}{2}+it\right)$ may be approximated by real zero polynomials
uniformly on compact sets, where $\xi$ denotes Riemann's xi-function. 

P\'olya-Schur's result generated a vast literature
on this subject and related topics, see \cite{BBCV} and references 
therein. Nevertheless, complete solutions to the fundamental problems of 
describing all linear operators preserving the set of 
real zero polynomials or, more generally, the set of polynomials with 
zero locus in a prescribed region $\Omega\subset \bC$, are yet to be 
found. Although many special cases and variations of these problems have 
been intensely studied for more than a century, to the best of our knowledge
they have been stated in the above general (and explicit) form only  recently by Craven-Csordas \cite{CC1} and Csordas \cite{cso}, see 
also \cite{ABH,Bo,BBCV,Br}.

This paper is part of a series \cite{Bo,BBS,BBS1,BBL,Br1,Br} devoted to these 
questions, their natural multivariate extensions and applications to geometric
function theory, matrix theory, probability theory, combinatorics, and 
statistical mechanics. Here we 
\begin{itemize}
\item[(1)] classify all linear partial differential operators in the 
$n$-th Weyl algebra that preserve stable, respectively real stable polynomials 
in $n$ variables;
\item[(2)] obtain a Lax type determinantal representation for linear operators 
in the first Weyl algebra ($n=1$) preserving real zero polynomials and a 
characterization in terms of their homogenized symbols;
\item[(3)] prove higher dimensional versions of P\'olya-Schur's theorem; 
\item[(4)] apply (1)--(3) to establish a Fischer-Fock duality for 
stability preservers, homotopic properties of their symbols and geometric 
interpretations extending P\'olya's ``curve theorem'' for all $n\ge 1$, 
stable multivariate generalizations of the Cauchy-Poincar\'e 
interlacing theorem, Schur-Mal\'o-Szeg\"o type convolution theorems in 
higher dimensions, and necessary and sufficient conditions for strict (real) 
stability preserving in one or several variables.
\end{itemize}

A nonzero univariate real polynomial with only real zeros is called 
{\em hyperbolic} while $f\in\CC[z]$ is called {\em stable} if 
$f(z) \neq 0$ for all $z \in \CC$ with $\Im(z) >0$. Hence a univariate 
real polynomial is stable if and only if it is hyperbolic. 
These classical concepts have several natural extensions to multivariate 
polynomials, see, e.g., four different definitions  in  \cite{KTM}. Below 
we concentrate on the most general notion:

\begin{definition}
A polynomial $f\in\bC[z_1,\ldots,z_n]$ is 
{\em stable} 
if  $f(z_1, \ldots, z_n) \neq 0$ for all $n$-tuples $(z_1, \ldots, z_n) \in \CC^n$ with $\Im(z_j) >0$, $1\le j\le n$.  If in addition $f$ has real coefficients it will be referred  to as {\em real stable}.
\end{definition}

Clearly, 
$f$ is stable (respectively, real stable) if and only if for all $\alpha \in \RR^n$ and $v \in \RR_+^n$ the 
univariate polynomial $f(\alpha + vt)\in\bC[t]$ is stable (respectively, hyperbolic), see Lemma \ref{lines} in \S \ref{s2}. In what follows we denote by $\HH_n(\CC)$, 
respectively $\HH_n(\RR)$, the 
set of stable, respectively real stable polynomials in $n$ variables.

Another fundamental extension of the notion of 
real-rootedness to higher dimensions stems from PDE theory. Namely, a   
homogeneous polynomial $p\in\bR[z_1,\ldots, z_n]$  is said to be {\em (G\aa rding) hyperbolic} with respect to a given 
vector $v \in \RR^n$ if $p(v) \neq 0$ and for all vectors $\alpha \in \RR^n$  the univariate polynomial $p(\alpha + vt)\in\bR[t]$ has only real zeros.  For background on (multivariate homogeneous) hyperbolic polynomials one may consult, 
e.g., \cite {ABG,Gaa,hor2}. In \S \ref{s52} we prove the following result describing the relation between real stable and hyperbolic polynomials. 

\begin{proposition}\label{pro-g}
Let $f\in\bR[z_1,\ldots, z_n]$ be of degree $d$ and let 
$f_{H}\in\bR[z_1,\ldots, z_{n+1}]$ be the (unique) homogeneous polynomial of degree $d$ such that 
$f_{H}(z_1, \ldots, z_n,1)=f(z_1,\ldots,z_n)$. Then $f\in\HH_n(\bR)$ 
if and only if 
$f_H$ is hyperbolic with respect to every vector $v\in\bR^{n+1}$ such that $v_{n+1}=0$ and $v_i>0$, $1 \leq i \leq n$. 
\end{proposition}
  
It is worth mentioning that real stable multivariate polynomials appear  already  in Theorem 1 of the foundational article   \cite {Gaa} by G\aa rding and that stable multivariate entire functions can be found in Chap.~IX of Levin's book \cite {Le}. 
   
Let $\A_n[\CC]$ be the Weyl algebra of all finite
order linear differential 
operators with polynomial coefficients on $\CC[z_1, \ldots, z_n]$. Recall the
standard multi-index notation $z^{\al}=z_1^{\al_1}\cdots z_{n}^{\al_n}$, 
$\partial^{\al}=\partial_1^{\al_1}\cdots \partial_n^{\al_n}$, where
$z=(z_1,\ldots,z_n)$, 
$\al=(\al_1,\ldots,\al_n)\in\bN^n$ and 
$\partial_i= {\partial}/{\partial z_i}$ for $1\le i\le n$.
Then each operator $T\in\A_n[\bC]$ may be (uniquely) represented as 
\begin{equation}\label{dop}
T= \sum_{\alpha,\beta \in \NN^n} a_{\alpha \beta} 
z^\alpha\partial^\beta,
\end{equation}
where $a_{\alpha \beta} \in \CC$ is nonzero only for a finite number of pairs $(\alpha, \beta)$. Let further 
$\A_n[\RR]$ be the set of all $T \in \A_n[\CC]$ with $a_{\alpha \beta} \in \RR$ for all 
$\alpha, \beta \in \NN^n$. 
A {\em nonzero} differential operator  $T \in \A_n[\CC]$ is called {\em stability preserving} 
if $T : \HH_n(\CC) \rightarrow \HH_n(\CC)\cup \{0\}$ and it is said to be 
{\em real stability preserving} if  $T : \HH_n(\RR) \rightarrow \HH_n(\RR)\cup \{0\}$.

Given $T$ of the form  \eqref{dop} define its  {\em symbol} $F_T(z,w)$ to be the polynomial in 
$\CC[z_1,\ldots,z_n,w_1,\ldots,w_n]$ given by $F_T(z,w)=\sum_{\alpha,\beta } 
a_{\alpha \beta} z^\alpha w^\beta$.

The first main results of this paper are the following characterizations of the multiplicative submonoids 
$\A_n(\CC) \subset \A_n[\CC]$ and  $\A_n(\RR) \subset \A_n[\RR]$ consisting of
all stability preservers and real stability preservers, respectively. 

\begin{theorem}\label{maincomplex}
Let $T \in \A_n[\CC]$. Then 
$
T \in \A_n(\CC)$ if and only if 
$F_T(z,-w) \in \HH_{2n}(\CC).
$
\end{theorem} 

 \begin{theorem}\label{mainreal}
Let $T \in \A_n[\RR]$. Then 
$
T \in \A_n(\RR)$ if and only if 
$F_T(z,-w) \in \HH_{2n}(\RR)$.
\end{theorem} 

It is interesting to note that Theorems~\ref{maincomplex} and \ref{mainreal} essentially assert that {\em finite order stability (respectively, real stability) preservers in $n$ variables are generated by stable (respectively, real stable) polynomials in $2n$ variables via the symbol map.} Geometric interpretations of these statements in terms of symbol surfaces are given in \S \ref{ss-polya}.

To prove the above theorems  we need to generalize  a large number of notions and results for univariate stable and hyperbolic polynomials to the multivariate case.    
Let $\alpha_1 \leq \alpha_2 \leq \cdots \leq \alpha_n$ and 
$\beta_1 \leq \beta_2 \leq \cdots \leq \beta_m$ be the zeros (counted with
multiplicities) of two 
given polynomials $f,g \in \HH_1(\RR)$ with $\deg f=n$ and  $\deg g=m$. We say that these zeros 
{\em interlace} if they can 
be ordered so that either   $\alpha_1 \leq \beta_1 \leq \alpha_2 \leq \beta_2 \leq 
\cdots$ or $\beta_1 \leq \alpha_1 \leq \beta_2 \leq \alpha_2 \leq 
\cdots$, in which case one clearly must have $|n-m|\le 1$. Note that by our convention, the zeros of any two polynomials 
of degree $0$ or $1$ interlace. It is not difficult to show that if   
the zeros of $f$ and $g$ interlace then 
the {\em Wronskian} 
$W[f,g]:=f'g-fg'$ is either nonnegative or nonpositive on the whole real axis 
$\RR$, see, e.g., \cite{RS}. In the case when  $W[f,g] \leq 0$ we say that $f$ and $g$ 
are  in  {\em proper position}, denoted $f \ll g$. 
For technical reasons we also say that the zeros of the polynomial $0$ 
interlace the zeros of any (nonzero) hyperbolic polynomial and write 
$0 \ll f$ and $f \ll 0$.
Note 
that if $f,g$ are (nonzero) hyperbolic polynomials such that $f \ll g$ and $g \ll f$ then $f$ and $g$ must be constant multiples 
of each other, that is, $W[f,g]\equiv 0$. 

The following theorem is a version of the classical Hermite-Biehler 
theorem \cite{RS}.

\begin{theorem}[(Hermite-Biehler theorem)]
Let $h := f +ig \in \CC[z]$, where $f,g \in \RR[z]$. Then 
$h \in \HH_1(\CC)$ if and only if $g \ll f$. 
\end{theorem}  

The Hermite-Biehler theorem gives an indication about how one should   
generalize the concept of 
interlacing  to higher dimensions: 

\begin{definition}\label{d-proper}
Two polynomials $f,g \in \RR[z_1,\ldots,z_n]$ are in 
{\em proper position}, denoted $f \ll g$, if $g+if \in \HH_n(\CC)$.
\end{definition} 

Equivalently, $f$ and $g$ are in proper 
position if and only if for all $\alpha \in  \RR^n$ and $v \in \RR_+^n$ the 
univariate polynomials 
$f(\alpha+vt), g(\alpha +vt)\in\bR[t]$ are in proper position.  
It also follows that  
$f, g \in \HH_n(\RR)\cup \{0\}$  whenever $f \ll g$, see Corollary \ref{hhC}
in \S \ref{s2}.

The next (also classical) result is often attributed to Obreschkoff 
\cite{obreschkoff} and sometimes referred to as the 
Hermite-Kakeya-Obreschkoff theorem \cite{RS}. 

\begin{theorem}[(Obreschkoff theorem)]\label{1Obreschkoff}
Let $f,g \in \RR[z]$. Then 
$\alpha f + \beta g \in \HH_1(\RR)\cup \{0\}$ for all $\alpha, \beta \in \RR$ 
if and only if either 
$f \ll g$, $g \ll f$ or  $f=g\equiv 0$.      
\end{theorem}

We extend this theorem to polynomials in several variables as follows.

\begin{theorem}\label{nObreschkoff}
Let $f,g \in \RR[z_1, \ldots, z_n]$. Then 
$\alpha f + \beta g \in \HH_n(\RR) \cup \{0\}$ for all 
$\alpha, \beta \in \RR$
if and only if either $f \ll g$, $g \ll f$ or $f=g\equiv 0$. 
\end{theorem}

Recall that an infinite  sequence of real numbers $\lambda : \NN \rightarrow \RR$ 
is called a {\em multiplier sequence (of the first kind)} if the associated
linear operator $T$ on $\bC[z]$ 
defined by $T(z^n)=\lambda(n) z^n$, for all $n\in\bN$, is a hyperbolicity preserver, i.e.,   
$T: \HH_1(\RR) \rightarrow \HH_1(\RR)\cup \{0\}$. Any linear 
operator $T:\CC[z_1,\ldots,z_n]\to \CC[z_1,\ldots,z_n]$ can be represented as a formal power series in $\partial$ with polynomial coefficients. Indeed, this may be
proved either by induction  or by 
invoking Peetre's abstract characterization of differential operators
\cite{Pe}. Note also that  in general a 
multiplier sequence is represented by  an infinite order 
differential operator with polynomial coefficients.

In \cite{PS} P\'olya and Schur gave the following characterization  
of multiplier sequences of the first kind. 

\begin{theorem}[(P\'olya-Schur theorem)]\label{ps}
Let $\lambda : \NN \rightarrow \RR$ be a sequence of real numbers and 
$T: \RR[z] \rightarrow \RR[z]$ be the corresponding (diagonal) linear operator.
 Define $\Phi(z)$ to be the formal power series 
$$
\Phi(z)= \sum_{k=0}^\infty \frac{\lambda(k)}{k!}z^k.
$$
The following assertions are equivalent:
\begin{itemize}
\item[(i)] $\lambda$ is a multiplier sequence,  
\item[(ii)] $\Phi(z)$ defines an entire 
function which is the limit, uniformly on compact sets, of 
polynomials with only real zeros of the same sign, 
\item[(iii)] Either $\Phi(z)$ or $\Phi(-z)$ is an entire function 
that can be written as
$$
C z^n e^{az} \prod_{k=1}^\infty (1+ \alpha_k z),
$$ 
where $n \in \NN$, $C \in \RR$, $a,\alpha_k \geq 0$ for all $k \in \NN$ and 
$\sum_{k=1}^\infty \alpha_k < \infty$, 
\item[(iv)] For all nonnegative integers $n$ the  polynomial 
$T[(1+z)^n]$ is hyperbolic with all zeros of the same sign. 
\end{itemize}
\end{theorem} 

We introduce a natural higher dimensional analog of the notion of 
multiplier sequence and completely characterize all multivariate multiplier 
sequences as well as those that can be represented as finite order differential
operators. For this we need the following notation. Given an integer $n\ge 1$ 
and $\alpha, \beta \in \NN^n$ we write $\alpha \leq \beta$ if 
$\alpha_i \leq \beta_i$ for $1 \leq i \leq n$. Let further
$|\alpha|= |\alpha_1|+\cdots+|\alpha_n|$, 
$\alpha^\beta = \alpha_1^{\beta_1}\cdots \alpha_n^{\beta_n}$, 
$\alpha! = \alpha_1! \cdots \alpha_n!$, $(\beta)_{\alpha}=0$ if 
$\alpha \nleq \beta$ and 
$(\beta)_{\alpha}=\beta! / (\beta-\alpha)!$ otherwise. 

\begin{definition}\label{d-mult-ms}
A function $\lambda : \NN^n \rightarrow \RR$ is a (multivariate) 
{\em multiplier sequence} if the corresponding (diagonal) 
linear operator $T$ defined 
by $T(z^{\alpha}) = \lambda(\alpha)z^{\alpha}$, for all $\al\in\bN^n$, is a real stability preserver, 
that is, 
$T : \HH_n(\RR) \rightarrow \HH_n(\RR)\cup \{0\}$. 
\end{definition}

The following theorem completely describes multivariate  
multiplier sequences.

\begin{theorem}\label{multichar}
Consider an arbitrary map $\lambda : \NN^n \rightarrow \RR$. Then $\lambda$ is a multivariate multiplier 
sequence if and only if there exist usual (univariate) 
multiplier sequences $\lambda_i : \NN \rightarrow \RR$, for all $1\le i\le n$, 
such that 
$$
\lambda(\alpha) = \lambda_1(\alpha_1)\lambda_2(\alpha_2)
\cdots \lambda_n(\alpha_n), \quad \mbox{ for all } \al=(\al_1,\ldots,\al_n)\in\bN^n,
$$
and either  $\lambda(\alpha)\lambda(\beta) \geq 0 $  for all $\alpha, \beta \in\bN^n$,  or $(-1)^{|\alpha|+|\beta|}\lambda(\alpha)\lambda(\beta)\geq 0$  for all $\alpha, \beta\in\bN^n$.
\end{theorem}  

Note that Theorem \ref{multichar} is a negative result, since it shows that the only multivariate multiplier sequences are those that one would expect exist: products of univariate ones. 

We next characterize all multiplier sequences that are 
finite order differential operators, i.e., those whose symbols are 
(finite degree) polynomials:  
 
\begin{theorem}\label{finitemulti}
 Given a map $\lambda : \NN^n \rightarrow \RR$,  let 
$T$ be the corresponding (diagonal) linear operator. Then $T \in \A_n(\RR)$ if 
and only if $T$ has a symbol $F_T(z,w)$ of the form
$$
F_T(z,w)=f_1(z_1w_1)f_2(z_2w_2)\cdots f_n(z_nw_n),
$$
where  $f_i(t)$ is a  polynomial with only  real and 
nonpositive zeros for all $1\le i\le n$.  
\end{theorem}

\begin{remark}
Note that Theorem~\ref{finitemulti} combined with well-known properties
of univariate multiplier sequences (cf.~Lemma~\ref{basicmult} below) implies
in particular that if $\lambda : \NN^n \rightarrow \RR$ is a finite order 
multivariate multiplier sequence, then there exists $\gamma \in \NN^n$ such 
that $\lambda(\alpha) = 0$ 
for $\alpha < \gamma$ and either $\lambda(\alpha)>0$ for 
all $\alpha \geq \gamma$ or $\lambda(\alpha)<0$
for all $\alpha \geq \gamma$. Note also that for $n=1$ 
Theorem~\ref{finitemulti} gives an alternative description of finite order
multiplier sequences that complements P\'olya-Schur's Theorem \ref{ps}.
\end{remark}

Our next result is a vast generalization of the following classical theorem
\cite{RS}.

\begin{theorem}[(Hermite-Poulain-Jensen theorem)]
Let $p(z) = \sum_{k=0}^{n}a_kz^k\in \RR[z]$ be nonzero and let
$T=\sum_{k=0}^{n}a_k{d^k}/{dz^k}\in 
\A_1[\RR]$. 
Then $T \in \A_1(\RR)$ if and only if $p \in \HH_1(\RR)$.
\end{theorem} 

The natural setting for our extension is the {\em Fischer-Fock space} 
$\FS_n$ \cite{Fi1,Fi2,Fo1,Fo2}, also called the 
{\em Bargmann-Segal space} \cite{Bargmann,Ba2,Se} or the 
{\em Newman-Shapiro space} \cite{NS1,NS2,NS3,H-S}, which is the Hilbert 
space of holomorphic functions $f$ on $\CC^n$ such that 
$$
\| f \|^2 = \sum_{\alpha\in \NN^n} \alpha!|a(\alpha)|^2 = 
\pi^{-n}\int |f(z)|^2 e^{-|z|^2} dz_1 \wedge \cdots \wedge dz_n <\infty.
$$ 
Here $\sum_\alpha a(\alpha)z^\alpha$ is the Taylor expansion of $f$. The inner 
product in $\FS_n$ is given by 
\begin{equation}\label{ff-p}
\langle f, g \rangle  = 
\pi^{-n}\int f(z)\overline{g(z)} e^{-|z|^2} dz_1 \wedge \cdots \wedge dz_n 
\end{equation}
and one can easily check  that monomials 
$\{ z^{\alpha}/\sqrt{\alpha!}\}_{\alpha \in \NN^n}$ 
form an  orthonormal basis. From this it follows that for $1 \leq i \leq n$
one has 
$$
\langle \partial_i z^\alpha, z^\beta \rangle = 
\alpha!\delta_{\alpha = \beta+e_i} =\langle z^\alpha, z_iz^\beta \rangle, $$
where $\delta$ is the Kronecker delta and $e_i$ denotes the $i$-th 
standard generator of the lattice $\bZ^n$. Hence, if 
$T = \sum_{\alpha, \beta}a_{\alpha \beta} z^\alpha\partial^\beta \in 
\A_n[\CC]$ then 
\begin{equation*}
 \langle T(f), g \rangle= \sum_{\alpha,\beta}a_{\alpha \beta}
\langle z^\alpha\partial^\beta f, g\rangle 
= \sum_{\alpha,\beta}a_{\alpha \beta}
\langle f, z^\beta\partial^\alpha g \rangle 
= \langle f, \sum_{\alpha,\beta}\overline{a_{\beta \alpha}}
z^\alpha\partial^\beta g \rangle. 
\end{equation*}
Therefore, the formal Fischer-Fock dual (or adjoint) operator of $T$ is 
given by 
$T^* =  \sum_{\alpha,\beta}\overline{a_{\beta \alpha}} 
z^\alpha\partial^\beta$, cf.~\cite{NS2}. Note that for $1\le i\le n$ 
the dual of $\pa_i$ is the operator given by multiplication with $z_i$ and 
that diagonal operators (in the standard monomial basis) are self-dual. In
particular, if $T$ is a multiplier sequence then $T^*=T$.

\begin{remark}
The Fischer-Fock space $\FS_n$ was used by Dirac to define second 
quantization \cite{dir} and its inner product has since been rediscovered in 
various contexts, e.g.~in number theory where the corresponding norm is 
known as the Bombieri norm \cite{BBEM,Rez}. Further important properties 
of $\FS_n$  
such as its (Bergman-Aronszajn) reproducing kernel and the 
Newman-Shapiro Isometry Theorem may be found in \cite{NS1,NS2,NS3}. 
We should also point out that in e.g.~$\calD$-module
theory and microlocal Fourier analysis \cite{Malgr} one usually works with 
the inner product on $\FS_n$ defined by 
$\langle f(z),g(z)\rangle_{d}=\langle f(iz), g(iz)\rangle$, where 
$\langle \cdot,\cdot\rangle$ is as in \eqref{ff-p}. Note that the dual 
operator of $\pa_i$ with respect to $\langle\cdot,\cdot\rangle_{d}$
is the operator given by multiplication with $-z_i$.
\end{remark}

In \S \ref{ss-polya} we give a geometric interpretation 
and proof of the fact that the 
duality map with respect to the above scalar product preserves both 
$\A_n(\CC)$ and $\A_n(\RR)$. More precisely, from 
Theorems~\ref{maincomplex}--\ref{mainreal} we
deduce the following natural property:

\begin{theorem}[(Duality theorem)]\label{adjoint}
Let $T \in \A_n[\CC]$. Then $T \in \A_n(\CC)$ if and only if 
$T^* \in \A_n(\CC)$. Similarly, if $T \in \A_n[\RR]$ then 
$T \in \A_n(\RR)$ if and only if $T^* \in \A_n(\RR)$.
\end{theorem}

We conclude this introduction with a series of examples of real stable 
polynomials and 
various applications of our results. Further interesting examples of 
multi-affine stable and real stable polynomials can be found in 
e.g.~\cite{BBS1,BBL,Br1,COSW}. 

\begin{proposition} \label{pencil}
 Let $A_1, \ldots, A_n$ be positive semidefinite $m \times m$ matrices 
and  let $B$ be a complex Hermitian  $m \times m$ matrix.  Then the polynomial 
\begin{equation}\label{determ}
f(z_1,\ldots,z_n)=\det\left(\sum_{i=1}^nz_iA_i+ B\right)
\end{equation}
is either real stable or identically zero.   
\end{proposition}

A proof of the above proposition is given in \S \ref{s-ap}. 
Using this result and Theorem~\ref{maincomplex} we obtain 
a multidimensional generalization of the Cauchy-Poincar\'e interlacing 
theorem: see Theorem~\ref{C-P} in \S \ref{s51}.

The Lax conjecture \cite{Lax} for (G\aa rding) hyperbolic polynomials in three 
variables has recently been settled by Lewis, Parillo and Ramana 
\cite{LPR}. Their proof relies on the results of Helton and 
Vinnikov \cite{HV}. Applications of these results to e.g.~hyperbolic 
programming and convex optimization may be found in \cite{Re}. 
In \S \ref{s52} we prove the following converse to 
Proposition~\ref{pencil} in the
case $n=2$ and thus establish a natural analog of the Lax conjecture for
real stable polynomials in two variables.  

\begin{theorem}\label{laxlike}
Any real stable polynomial in two variables $x, y$ can be written as 
$\pm\det(xA+yB+C)$ where $A$ and $B$ are positive semidefinite matrices 
and $C$ is a symmetric matrix of the same order. 
\end{theorem}  

\begin{remark}
A characterization of real stable polynomials in an arbitrary number of 
variables has recently been 
obtained in \cite{Br1}.
\end{remark}

Combining Theorem~\ref{laxlike} with Theorem~\ref{mainreal} we get two new
descriptions of finite order linear preservers of hyperbolicity (i.e., 
univariate real stability), namely a determinantal
characterization and one in terms of homogenized operator symbols: see 
Theorems~\ref{crit1} and~\ref{crit2} in \S \ref{s53}. 

Further applications of our results include multivariate Schur-Mal\'o-Szeg\H{o} 
composition formulas and closure properties under the 
Weyl product of (real) stable polynomials 
(\S \ref{s54}), a unified treatment of P\'olya type ``curve theorems'' as well
as multivariate extensions (\S \ref{ss-polya}), and necessary and 
sufficient 
criteria for strict stability and strict real stability preservers 
(\S \ref{s-strict}). 


\subsubsection*{Brief excursion around the literature.} 
The study of univariate stable polynomials was initiated by Hermite in the 
1860's and continued by Laguerre, Maxwell, Routh, Hurwitz and many others in 
the second half of the XIX-th century. The contributions of the classical 
period are well summarized in \cite{Gant,PSz,RS}. Important results on 
stability of entire functions were obtained in the mid XX-th century by 
e.g.~Krein, 
Pontryagin, Chebotarev, Levin \cite{Le}. Modern achievements  in this area 
can be found in \cite {SaOl} and references therein. Much less seems to be 
known concerning multidimensional stability. In control theory one can name 
a series of papers by Kharitonov {\em et al.}
\cite {KTM} with numerous references to the earlier literature on this 
topic.  Another origin of interest to multivariate stable polynomials comes 
from an unexpected direction, namely the Lee-Yang theorem on ferromagnetic 
Ising models, the Heilmann-Lieb theorem for monomer-dimer systems and their
various generalizations \cite{HL,LY,LS}. Combinatorial theory provides yet 
another rich source of stable polynomials as multivariate spanning tree 
polynomials and generating 
polynomials for various classes of matroids turn out to be stable (cf., e.g., 
\cite{Br1,COSW}). Multivariate stable polynomials were recently used in 
\cite{Gu2} to generalize  and reprove in a unified manner a number of   
classical  conjectures, including  the van der Waerden and Schrijver-Valiant 
conjectures, and in \cite{BBS1} to solve some long-standing  
conjectures of Johnson and Bapat in matrix theory. Further recent  
contributions include \cite{BBS}, where a complete classification of 
linear preservers of univariate
polynomials with only real zeros -- and, more generally, of univariate
polynomials with zeros only in a closed circular
domain or on the boundary of such a domain -- has been obtained, thus 
solving an
old open problem going back to Laguerre \cite{lag}  and
P\'olya-Schur \cite{PS}. Let us finally mention that
real stable polynomials have also found remarkable applications in 
probability theory
and interacting particle systems. Indeed, these polynomials were recently 
used in \cite{BBL} to develop a theory of negative dependence for the class 
of strongly Rayleigh probability measures, which contains several important 
examples such as
uniform random spanning tree measures, fermionic/determinantal measures, 
balls-and-bins measures and distributions for symmetric exclusion processes. \\[2ex]
{\bf Acknowledgements}. 
It is a pleasure to thank the American Institute of Mathematics for
hosting the ``P\'olya-Schur-Lax Workshop'' on these themes in Spring 2007. We would also like to thank an anonymous referee for several suggestions on improving the exposition of this paper.

\section{Basic properties and generalized Hermite-Kakeya-Obreschkoff 
Theorem}\label{s2}

The following criterion for (real) stability is an easy consequence of the 
definitions.

\begin{lemma}\label{lines}
Let $f \in \CC[z_1,\ldots,z_n]$. Then 
$f \in \HH_n(\CC)$ if and only if 
$f(\alpha + vt) \in \HH_1(\CC)$ 
for all $\alpha \in \RR^n$ and $v \in \RR_+^n$. 
\end{lemma}

The next lemma extends the Hermite-Biehler theorem to the multivariate case 
and provides a useful alternative description of the proper 
position/``interlacing'' property for multivariate polynomials.

\begin{lemma}\label{addz}
Let $f, g \in \RR[z_1, \ldots, z_n]$ and let $z_{n+1}$ be a new indeterminate.
Then $f \ll g$ if and only 
if $g+z_{n+1}f \in \HH_{n+1}(\RR)$. Moreover, if $f \in \HH_n(\RR)$ then 
$f \ll g$ if and only if 
$$
\Im\left(\frac {g(z)}{f(z)}\right) \geq 0
$$ 
whenever $\Im(z_i) >0$ for all $1\leq i \leq n$.
\end{lemma}

\begin{proof}
The ``if'' direction is obvious. Suppose that $f \ll g$ and that 
$z_{n+1}=a+ib$, where $a \in \RR$ and $b \in \RR_+$. Then 
by Lemma \ref{lines} we have that $f(\alpha + vt) \ll g(\alpha + vt)$ for 
all $\alpha \in \RR^n$ and $v \in \RR_+^n$. By Obreschkoff's theorem the 
zeros of $g(\alpha + vt)+af(\alpha + vt)$ and $bf(\alpha + vt)$ 
interlace (both cannot be identically zero). Moreover, 
$$W(bf(\alpha + vt),g(\alpha + vt)+af(\alpha + vt))= 
bW(f(\alpha + vt)),g(\alpha + vt)).
$$ Thus   
$bf(\alpha + vt)\ll g(\alpha + vt)+af(\alpha + vt)$ for all 
$\alpha \in \RR^n$ and $v \in \RR_+^n$, which by Lemma \ref{lines} implies that  
$g+(a+ib)f \in \HH_{n}(\CC)$. But 
$g+z_{n+1}f$ clearly has real coefficients so 
$g+z_{n+1}f \in \HH_{n+1}(\RR)$. The final statement of the lemma is 
a simple consequence of the above arguments.
\end{proof}

\begin{lemma}\label{althurwitz}
Suppose that $f_j \in \CC[z_1,\ldots, z_n]$ for all $j\in\bN$  are 
nonvanishing 
in an open set $U \subseteq \CC^n$ and that $f$ is the limit, uniformly on 
compact sets, of the sequence $\{f_j\}_{j\in\bN}$. Then $f$ is either nonvanishing 
in $U$ or it is identically 
equal to $0$.
\end{lemma}

\begin{proof}
The lemma follows from the multivariate version of Hurwitz' theorem on the 
continuity of zeros of analytic functions, see, e.g., \cite{COSW}.
\end{proof}

Let $f(z_1, \ldots, z_n) \in \HH_n(\CC)$, $\alpha \in \RR$ and 
$\lambda >0$.  Then 
$f(\alpha +\lambda z_1, \ldots, z_n) \in \HH_n(\CC)$. By letting 
$\lambda\to 0$ we have by Lemma \ref{althurwitz} that 
$f(\alpha,z_2, \ldots, z_n) \in \HH_{n-1}(\CC)\cup \{0\}$.  

\begin{corollary}\label{hhC}
For each $n\in\bN$ one has
$$
\HH_n(\CC) = \{ g+if : f,g \in \HH_n(\RR)\cup \{0\}, f \ll g \}. 
$$ 
\end{corollary}
\begin{proof} 
The only novel part is that $f,g \in \HH_n(\RR)\cup \{0\}$ 
whenever $f \ll g$. This follows from Lemma \ref{addz} and Lemma 
\ref{althurwitz} when we let $z_{n+1}$ tend to $0$ and $\infty$, 
respectively. 
\end{proof}

We are ready to prove our multivariate Obreschkoff theorem.  

\begin{proof}[of Theorem \ref{nObreschkoff}]
Suppose that $f \ll g$. By Corollary \ref{hhC} we have 
$g \in \HH_n(\RR) \cup \{0\}$ so we can normalize and set $\beta=1$. By 
Lemma \ref{addz} we have 
$g + z_{n+1}f \in \HH_{n+1}(\RR) \subset \HH_{n+1}(\CC)$, so by letting 
$z_{n+1}=i + \alpha$ with $\alpha \in \RR$ we have 
$g + \alpha f +if \in \HH_{n}(\CC)$, i.e., 
$f \ll g + \alpha f$. From Corollary \ref{hhC} again it follows that 
$g + \alpha f \in \HH_n(\RR) \cup \{0\}$, as was to be shown. 

To prove the converse statement suppose that we do not have $f=g\equiv 0$.  
If $\alpha f + \beta g \in \HH_n(\RR) \cup \{0\}$ 
for all $\alpha, \beta \in \RR$ then (by Lemma \ref{lines} and 
Obreschkoff's theorem) for all 
$\gamma \in \RR^n$ and $v \in \RR_+^n$ we have either 
$f(\gamma + vt) \ll g(\gamma + vt)$ or $f(\gamma + vt) \gg g(\gamma + vt)$. 
If both instances  occur for different vectors, i.e., 
$f(\gamma_1 + v_1t) \ll g(\gamma_1 + v_1t)$ and 
$f(\gamma_2 + v_2t) \gg g(\gamma_2 + v_2t)$ for some $\ga_1,\ga_2\in\bR^n$
and $v_1,v_2\in\RR_+^n$, then by continuity arguments there exists 
$\tau\in[0,1]$ such that 
$f(\gamma_\tau + v_\tau t) \ll g(\gamma_\tau + v_\tau t)$ and 
$f(\gamma_\tau + v_\tau t) \gg g(\gamma_\tau + v_\tau t)$, where  
$\gamma_\tau:=\tau\ga_1+(1-\tau)\ga_2 \in \RR^n$ and 
$v_\tau:=\tau v_1+(1-\tau)v_2 \in \RR_+^n$ . 
This means that $f(\gamma_\tau + v_\tau t)$ and $g(\gamma_\tau  + v_\tau t)$ 
are constant 
multiples of each other, say $f(\gamma_\tau + v_\tau t)
= \lambda g(\gamma_\tau + v_\tau t)$ for 
some $\lambda \in \RR$. By hypothesis we have  
$h:=f-\lambda g \in \HH_n(\RR) \cup \{0\}$ and $h(\gamma_\tau +v_\tau t) 
\equiv 0$, in particular $h(\gamma_\tau +iv_\tau)=0$. 
Since $v_\tau\in\RR_+^n$ it follows that 
$h \equiv 0$ and $f= \lambda g$. 
Consequently,  
if both instances occur we  have $f \ll g$ for trivial reasons. Thus we may assume that 
only one of them occurs.  But then the conclusion follows from Lemma \ref{lines}.  
\end{proof}

Define 
$$
\HH_n(\CC)^{-}= \{ f \in \CC[z_1,\ldots,z_n]: f(z_1, ,\ldots,z_n) \neq 0 
\mbox{ if } \Im(z_i) < 0 \mbox{ for } 1 \leq i \leq n\}. 
$$
Clearly, $f(z) \in  \HH_n(\CC)$ if and only if $f(-z) \in  \HH_n(\CC)^-$. 
Hence by Corollary \ref{hhC} and 
Lemma \ref{lines}  we have $h:=g+if \in \HH_n(\CC)^-$ with 
$f,g \in \RR[z_1,\ldots,z_n]$ if and only 
if $g \ll f$. 

\begin{proposition}
For any $n\in\bN$ the following holds:
$$\HH_n(\CC) \cap \HH_n(\CC)^{-} = \CC \HH_n(\RR) 
         := \{ cf: c \in \CC, f \in \HH_n(\RR) \}.$$
\end{proposition}
\begin{proof}
Suppose that $h = g+if \in \HH_n(\CC) \cap \HH_n(\CC)^{-}$. 
By Corollary \ref{hhC} we have $f \ll g$ and $g \ll f$. Hence  
for all $\alpha \in \RR^n$ and $v \in \RR_+^n$ we also have 
$f(\alpha + vt) \ll g(\alpha + vt)$ and $g(\alpha + vt) \ll f(\alpha + vt)$. 
This means that $f(\alpha + vt)$ and $g(\alpha + vt)$ are constant 
multiples of each other, say $f(\alpha + vt) = \lambda g(\alpha + vt)$. 
By the multivariate Obreschkoff theorem we have 
that $f-\lambda g \in \HH_n(\RR)\cup \{0\}$. Since 
$(f-\lambda g)(\alpha + vi)=0$ we must have $f-\lambda g\equiv 0$, i.e., 
$h = (1+i\lambda)g \in \CC \HH_n(\RR)$. 
\end{proof}  

\begin{lemma}\label{cones}
Let $f \in \HH_n(\RR)$. Then the sets 
$$
\{g \in \HH_n(\RR): f \ll g\} \mbox{ and } 
\{g \in \HH_n(\RR): f \gg g\}
$$
are nonnegative cones, i.e., they are closed under nonnegative 
linear combinations.
\end{lemma}

\begin{proof}
Let $f \in \HH_n(\RR)$ and suppose that $f \ll g$ and $f \ll h$. Then by 
Lemma~\ref{addz} we have $\Im(g(z)/f(z)) \geq 0$ and 
$\Im(h(z)/f(z)) \geq 0$ whenever $\Im(z) > 0$. Hence if 
$\lambda, \mu \geq 0$ and $\Im(z) > 0$ then 
$\Im\left( (\lambda g(z) + \mu h(z))/f(z) \right) \geq 0$ and 
Lemma~\ref{addz} yields 
$f \ll \lambda g + \mu h$. The other assertion follows similarly. 
\end{proof}

\section{Classifications of multivariate multiplier sequences and finite
order ones}

\subsection{Univariate and multivariate multiplier sequences}

Let us first recall a few well-known properties of (usual) univariate 
multiplier sequences, see, e.g., \cite{CC1}. 

\begin{lemma}\label{basicmult}
Let $\lambda:\bN\to \bR$  
be a multiplier sequence. If $0 \leq i \leq j \leq k $ are such that 
$\lambda(i)\lambda(k) \neq 0$ then $\lambda(j) \neq 0$. Furthermore, either 
\begin{itemize} 
\item[(i)] all nonzero $\lambda(i)$ have the same sign, or
\item[(ii)] all nonzero entries of the sequence 
$\{(-1)^i\lambda(i)\}_{i\geq 0}$ have the same sign. 
\end{itemize}
\end{lemma}

In what follows we denote the standard basis in $\bR^n$ by 
$\{e_k: 1 \leq k \leq n\}$.

\begin{remark}\label{proj} 
Suppose that $\alpha \in \NN^n$, $1 \leq k \leq n$, 
$f(z_k) := \sum_{i=0}^N a_iz_k^i \in \HH_1(\RR)$ and assume that $\lambda$ is 
a multivariate multiplier sequence. Then 
$$
T(z^\alpha f(z_k)) = z^\alpha\sum_{i = 0}^N \lambda(\alpha+ie_k)a_iz_k^i \in \HH_n(\RR). 
$$
Hence the function $i \mapsto \lambda(\al+ie_k)$ is a univariate multiplier 
sequence. 
\end{remark} 

The proofs of our characterizations of multivariate multiplier sequences and
those of finite order build on a series of statements that we proceed to
describe.

\begin{lemma}\label{det}
Let $f(z_1,z_2) = a_{00} + a_{01}z_2 + a_{10}z_1 + a_{11}z_1z_2 \in 
\RR[z_1,z_2]\setminus \{0\}$. Then $f \in \HH_2(\RR)$ if and only if 
$\det (a_{ij}) \leq 0$. 
\end{lemma} 
\begin{proof}
Let $\alpha \in \RR$ and denote  by $A=(a_{ij})$ the matrix of coefficients of 
$f(z_1, z_2)$. 
Clearly, $f(z_1,z_2) \in \HH_2(\RR)$ if and only if 
$f(z_1+\alpha,z_2) \in \HH_2(\RR)$ and  
$f(z_1,z_2+\alpha) \in \HH_2(\RR)$. 
 We get the matrix corresponding to 
$f(z_1+\alpha,z_2)$ by adding  $\alpha$ times  the last row of $A$ 
to the first row of $A$, and we get the matrix corresponding to 
$f(z_1,z_2+\alpha)$ by adding   $\alpha$ times the last column of $A$ 
to the first column of $A$. Since the determinant is preserved   
under such  row and column operations we can assume that 
$A$ has one of the following forms: 
$$
\begin{pmatrix}
a_{00} & 0 \\
0 & a_{11}
\end{pmatrix}
, \ \
\begin{pmatrix}
0 & a_{01} \\
a_{10} & 0
\end{pmatrix}
, \ \
\begin{pmatrix}
a_{00} & a_{01} \\
0 & 0
\end{pmatrix}.
$$
Obviously, these matrices correspond to a polynomial $f(z_1,z_2) \in \HH_2(\RR)$ if and only 
if $\det (a_{ij}) \leq 0$.
\end{proof}

\begin{lemma}\label{det0}
Let $\lambda :\NN^n \rightarrow \RR$, $n \geq 2$, be a multivariate multiplier 
sequence and let $\gamma \in \NN^n$ and $1\leq i<j \leq n$. Then 
$$
\lambda(\gamma) \lambda(\gamma + e_i + e_j) = 
\lambda(\gamma+e_i)\lambda(\gamma + e_j).
$$
\end{lemma}

\begin{proof}
Without loss of generality we may  assume that $i=1$ and $j=2$. Let  
$f(z)=z^\gamma g(z) \in \HH_n(\RR)$, where 
$g(z_1,z_2)= a_{00} + a_{01}z_2 + a_{10}z_1 + a_{11}z_1z_2 \in \HH_2(\RR)$. It follows 
that the polynomial
$$
\lambda(\gamma)a_{00} + \lambda(\gamma + e_2)a_{01}z_2 + 
\lambda(\gamma+e_1)a_{10}z_1 
+ \lambda(\gamma + e_1 + e_2)a_{11}z_1z_2 
$$
is stable or identically zero. 
By choosing $A$ as 
$$
\begin{pmatrix}
1 & 1 \\
1 & 1
\end{pmatrix}
\mbox{ and } 
\begin{pmatrix}
1 & -1 \\
1 & -1
\end{pmatrix},
$$
respectively, we get by Lemma \ref{det} that 
$
\lambda(\gamma) \lambda(\gamma + e_1 + e_2) \leq  
\lambda(\gamma+e_1)\lambda(\gamma + e_2)
$
and 
$
\lambda(\gamma) \lambda(\gamma + e_1 + e_2) \geq  
\lambda(\gamma+e_1)\lambda(\gamma + e_2)
$, respectively, 
which proves the lemma.
\end{proof}

Given $\al,\be\in\bN^n$ with $\al\le \be$ set $[\al,\be]:=\{\ga\in\bN^n:
\alpha \leq \gamma \leq \beta\}$.

\begin{lemma}\label{interval}
Let $\lambda :\NN^n \rightarrow \RR$ be a multivariate multiplier sequence. 
If $\al,\be\in\bN^n$ are such that  $\lambda(\alpha)\lambda(\beta) \neq 0$
and $\ga\in[\al,\be]$ then $\lambda(\gamma) \neq 0$. 
\end{lemma}

\begin{proof}
We use  induction on  
$\ell=|\beta|-|\alpha|$, the length of the interval $[\alpha,\beta]$. 
The cases 
$\ell = 0$ and $\ell=1$ are clear. By Remark \ref{proj} and Lemma \ref{basicmult} 
the result is true in the univariate case.  So  
we may assume that $\alpha$ and $\beta$ differ in more than one coordinate, 
i.e.,  that $\alpha + e_1, \alpha + e_2 \in [\alpha, \beta]$. 

If there exists an atom $\alpha + e_i \in [\alpha, \beta]$ such that 
$\lambda(\alpha + e_i) \neq 0$ then by induction we have that 
$\lambda$ is nonzero in $[\alpha+e_i, \beta]$. If $\alpha + e_j$ is 
another atom then $\alpha + e_i +e_j \in [\alpha+e_i, \beta]$ so 
$\lambda(\alpha + e_i +e_j) \neq 0$. Lemma \ref{det0} then gives that 
$\lambda(\alpha + e_j) \neq 0$. Thus, by induction, $\lambda$ is nonzero in 
$[\alpha+e_j, \beta]$ for all $\alpha + e_j \in [\alpha, \beta]$ and 
we are done.

In order to get a contradiction we may assume by the above that 
$\lambda(\alpha + e_i) =0$ for all 
$\alpha + e_i \in [\alpha, \beta]$.  Let $\gamma \in (\alpha, \beta]$ 
be a minimal element such that $\lambda(\gamma) \neq 0$. If $T$ is the 
(diagonal) linear operator associated to $\lambda$ then
$$
T(z^\alpha(1+z)^{\gamma - \alpha}) = 
\lambda(\alpha)z^\alpha + \lambda(\gamma)z^\gamma \in \HH_n(\RR). 
$$
By Lemma \ref{det0} we have that 
$\lambda(\alpha + e_i+ e_j) =0$ for all atoms $\alpha + e_i,  \alpha + e_j \in [\alpha,\beta]$. By Remark \ref{proj} and Lemma \ref{basicmult} we 
also have $\lambda(\alpha + me_i)=0$ for all $m \geq 1$ and atoms $\alpha + e_i \in [\alpha, \beta]$. It follows that 
$|\gamma| - |\alpha| \geq 3$. Now, if we set $z_i = t$ for all $i$ then 
by the above we obtain that the polynomial 
$$
\lambda(\alpha)t^{|\alpha|} + \lambda(\gamma)t^{|\gamma|}
$$
is hyperbolic in $t$, which is a contradiction since 
$|\gamma| - |\alpha| \geq 3$ and $\lambda(\alpha)\lambda(\gamma) \neq 0$.  
\end{proof} 

\begin{lemma}\label{formel}
Let $\lambda :\NN^n \rightarrow \RR$ be a multivariate multiplier sequence and 
suppose that $\lambda(\alpha) \neq 0$. Then 
\begin{equation}\label{formeln}
\frac{\lambda(\beta)}{\lambda(\al)}=\prod_{i=1}^n\frac{\lambda(\al+(\beta_i-\al_i)e_i)}{\lambda(\al)}
\end{equation}
for all $\beta \geq \alpha$. 
\end{lemma}

\begin{proof}
We claim that $
\lambda(\gamma) \lambda(\gamma + ae_i + be_j) = 
\lambda(\gamma+ae_i)\lambda(\gamma + be_j)
$ for all $\gamma \in \NN^n$ and positive integers $a$ and $b$. This follows easily  by induction on $a+b$ using Lemma~\ref{det0} and  Lemma \ref{interval}. We proceed to prove the \eqref{formeln} by induction on the number $r$ of nonzero entries of $\delta=\beta-\al$. The basis of induction, $r\leq 1$, is trivial. For $r\geq 2$ let $j$ be an index such that  $\delta_j>0$ and consider $\gamma=\al+\delta_je_j$. By Lemma \ref{interval} $\lambda(\gamma) \neq 0$. By induction 
$$
\frac{\lambda(\beta)}{\lambda(\al)}= \frac{\lambda(\beta)}{\lambda(\gamma)} \cdot \frac{\lambda(\gamma)}{\lambda(\al)} = \prod_{i=1}^n\frac{\lambda(\al+(\beta_i-\gamma_i)e_i)}{\lambda(\gamma)}\cdot \frac{\lambda(\al+(\gamma_i-\al_i)e_i)}{\lambda(\al)}.
$$
For $i\neq j$ the corresponding factor is 
$$
\frac{\lambda(\al+\delta_ie_i+\delta_j e_j)}{\lambda(\al+\delta_je_j)}\cdot \frac{\lambda(\al)}{\lambda(\al)}= \frac{\lambda(\al +\delta_ie_i)}{\delta(\al)}
$$
by the claim above. For $i=j$ the corresponding factor is 
$$
\frac{\lambda(\al+\delta_j e_j)}{\lambda(\al+\delta_je_j)}\cdot \frac{\lambda(\al+\delta_je_j)}{\lambda(\al)}= \frac{\lambda(\al +\delta_ie_j)}{\delta(\al)}
$$
and the lemma follows by induction. 

\end{proof}

If $\al=(\al_1,\ldots,\al_n)\in\bN^n$ and $\be=(\be_1,\ldots,\be_n)\in\bN^n$ 
we define two new vectors $\al\vee\be,\al\wedge\be\in\bN^n$ by setting
$\al\vee\be=(\max(\al_1,\be_1),\ldots,\max(\al_n,\be_n))$ and 
$\al\wedge \be=(\min(\al_1,\be_1),\ldots,\min(\al_n,\be_n))$.

\begin{lemma}\label{hathut}
Let $\lambda :\NN^n \rightarrow \RR$ be a multivariate multiplier sequence 
and suppose that $\lambda(\alpha)\lambda(\beta) \neq 0$. Then 
$\lambda(\alpha \vee \beta) \neq 0$ if and only if 
$\lambda(\alpha \wedge \beta) \neq 0$.  
\end{lemma}

\begin{proof}
If $\lambda(\alpha)\lambda(\beta) \neq 0$ and 
$\lambda(\alpha \wedge \beta) \neq 0$ then Lemma \ref{formel} and Lemma 
\ref{basicmult} imply  
that $\lambda(\alpha \vee \beta) \neq 0$. 
Suppose that 
$\lambda(\alpha)\lambda(\beta) \neq 0$ and 
$\lambda(\alpha \vee \beta) \neq 0$. We prove that 
$\lambda(\alpha \wedge \beta) \neq 0$ by induction on $|\alpha -\beta|$. 
If $|\alpha -\beta|=0$ there is nothing to prove. Also, if  
$\alpha$ and $\beta$ are comparable there is nothing to prove, 
so we may assume that there are indices $i$ and $j$ such that 
$\alpha_i < \beta_i$ and $\beta_j < \alpha_j$. Since 
$\alpha < \alpha + e_i \leq \alpha \vee \beta$ and 
$\beta < \beta + e_j \leq \alpha \vee \beta$ we have by Lemma 
\ref{interval} that $\lambda(\alpha + e_i)\lambda(\beta + e_j) \neq 0$. 
Consider the pairs $(\alpha + e_i, \beta+e_j)$,   
$(\alpha + e_i, \beta)$ and $(\alpha, \beta+e_j)$. The distance between each of them is smaller 
than $|\alpha-\beta|$, they all have to join $\alpha \vee \beta$, and the meets are 
$\alpha \wedge \beta +e_i+e_j$, $\alpha \wedge \beta +e_j$ and 
$\alpha \wedge \beta +e_i$ respectively, see Fig.~\ref{figge}. By induction we have that 
$\lambda(\alpha \wedge \beta +e_i)\lambda(\alpha \wedge \beta +e_j)
\lambda(\alpha \wedge \beta +e_i+e_j)\neq 0$. By Lemma \ref{det0} this 
gives 
$$
\lambda(\alpha \wedge \beta) = \frac{\lambda(\alpha \wedge \beta +e_i)
\lambda(\alpha \wedge \beta +e_j)}{\lambda(\alpha \wedge \beta +e_i+e_j)} 
\neq 0, 
$$
which is the desired conclusion.
\begin{figure}
$$
\vcenter{\xymatrix@R=10pt@C=10pt{
&&&&\alpha \vee \beta&&&&\\
&&&\alpha + e_i \de[ru]&&\de[lu] \beta +e_j&&&\\
&&\alpha \e[ru]&&\alpha \wedge \beta+e_i+e_j  \de[lu]\de[ru]&&\e[lu] \beta&& \\
&&&\alpha \wedge \beta + e_j\de[lu]\e[ru]&&\e[lu]\de[ru]\alpha \wedge \beta + e_i&&& \\
&&&&\alpha \wedge \beta \e[lu]\e[ru]&&&&
    }}
$$
\caption{\label{figge} Illustration of the induction step in  Lemma \ref{hathut}.}
\end{figure}
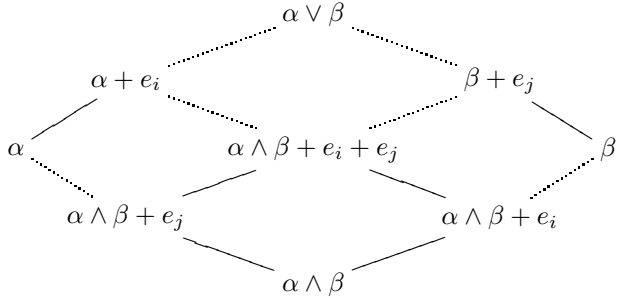
\end{proof}

Recall that the {\em support}
of a map $\lambda :\NN^n \rightarrow \RR$ is the set $\{\al\in\bN^n:\la(\al)\neq 0\}$.

\begin{lemma}
Let $\lambda :\NN^n \rightarrow \RR$ be a multivariate multiplier sequence. 
Then there exist univariate  multiplier sequences 
$\lambda_i:\NN\rightarrow \RR$, 
$1 \leq i \leq n$, such that 
$$
\lambda(\alpha)= \lambda_1(\alpha_1)\lambda_2(\alpha_2)
\cdots \lambda_n(\alpha_n),\quad \al=(\al_1,\ldots,\al_n)\in\bN^n.
$$
\end{lemma}
\begin{proof}
Let $S$ be the support of $\lambda$. By Lemma \ref{formel} and Remark \ref{proj} it suffices to prove that  $S$ has a unique minimal element. So far, by Lemma \ref{interval} and 
Lemma \ref{hathut} we know that $S$  is a disjoint union 
$S = \cup_{i=1}^m B_i$ of boxes $B_i = 
I_1^i\times \cdots \times  I_n^i$, 
where $I_n^i$ is an interval (possibly infinite) of nonnegative 
integers. Also, points in different boxes are incomparable.     

Suppose that $S$ does not have a unique minimal element. 
We claim that there exists an interval $[\alpha, \beta]$ such that 
$[\alpha, \beta]\cap S =\{\delta, \gamma\}$, where $\delta$ and $\gamma$ are 
in different boxes. We postpone the proof of this statement for a while 
and show first how it leads to a contradiction. Let $T$ be the (diagonal) 
linear 
operator associated to $\lambda$. We then have 
$$
T(z^\alpha(1+z)^{\beta-\alpha})= \lambda(\delta)z^\delta + 
\lambda(\gamma)z^\gamma.  
$$
Now $|\delta -\gamma| \geq 3$ since otherwise $\delta$ and $\gamma$ would 
be comparable or we would have 
$\gamma = \delta \wedge \gamma + e_i$ and 
$\delta = \delta \wedge \gamma + e_j$ for some $i$ and $j$. This is 
impossible by Lemma \ref{det0} since $\gamma$ and $\delta$ would then be 
in the same box. By assumption  we have that 
$$
\lambda(\delta)z^{\delta-\delta \wedge \gamma}+ 
\lambda(\gamma)z^{\gamma-\delta \wedge \gamma} \in \HH_n(\RR)
$$
so by setting all the variables in $z^{\delta-\delta \wedge \gamma}$  
equal to $t$ and setting all the variables in 
$z^{\gamma-\delta \wedge \gamma}$  
equal to $-t^{-1}$ (which we may since $z^{\delta-\delta \wedge \gamma}$ and
$z^{\gamma-\delta \wedge \gamma}$ contain no common variables) we obtain that 
$$
\lambda(\delta)t^{|\delta-\gamma|} \pm  
\lambda(\gamma)\in \HH_1(\RR). 
$$
This is a contradiction since $|\delta -\gamma| \geq 3$ and 
$\delta,\gamma\in S$ so 
$\la(\delta)\la(\ga)\neq 0$.

It remains to prove the claim. 
Let $d$ be the minimal distance between different boxes and 
suppose that $\delta$ and $\gamma$ are two points that realize the minimal 
distance. If $\kappa \in  
[\delta \wedge \gamma, \delta \vee \gamma]$ then 
$|\kappa -\delta| \leq d$ with equality only if $\kappa =\gamma$ and 
$|\kappa -\gamma| \leq d$ with equality only if $\kappa =\delta$. 
It follows that 
$[\delta \wedge \gamma, \delta \vee \gamma]\cap S =\{\delta, \gamma\}$.
\end{proof}

\subsection{Affine differential contractions and multivariate 
compositions}\label{s32}

For the proof of Theorem \ref{multichar} we need to establish first 
Theorem \ref{multischur} below, which is the main purpose of this section.

The proof of Theorem \ref{multischur} relies on some of the results obtained 
in \cite{LS}. Let 
us introduce the following notation. Given $a,b\in\bC$, $1\le i<j\le n$ and
\begin{equation*}
F(z_1,\ldots,z_n)=
\sum_{\al \in \NN^n} a_{\al}z^{\al}
\in\bC[z_1,\ldots,z_n]
\end{equation*}
let
\begin{equation}\label{not-contr}
F\!\left(z_1,\ldots,z_{i-1},az_i+b\frac{\partial}{\partial z_j},z_{i+1},
\ldots,z_j,\ldots,z_n\right)
\end{equation}
denote the polynomial
$$\sum_{\al \in \NN^n} a_{\al}z_1^{\al_1}\cdots 
z_{i-1}^{\al_{i-1}}\left(az_i+b\frac{\partial}{\partial z_j}\right)^{\al_i}
\!z_{i+1}^{\al_{i+1}}\cdots z_j^{\al_j}
\cdots z_n^{\al_n}.$$  

The next lemma follows from \cite[Lemma 2.3]{LS} by a rotation of the 
variables. 

\begin{lemma}[(\cite{LS})]\label{l-LS}
If $P_0(v),P_1(v)\in\bC[v]$ with $P_0(v)+xP_1(v)\neq 0$ for
$\Im (v)\ge c$ and $\Im (x)\ge d$ then
$$P_0(v)+\left(x-\frac{\partial}{\partial v}\right)P_1(v)\neq 0$$
for $\Im (v)\ge c$ and $\Im (x)\ge d$.
\end{lemma}

Using Lemma~\ref{l-LS} and the Grace-Walsh-Szeg\H{o} Coincidence Theorem 
\cite{RS} one can argue as in the proof of \cite[Proposition 2.2]{LS} to show:

\begin{proposition}[(\cite{LS})]\label{g-LS}
Let $(c_1,\ldots,c_n)\in\bR^n$ and $F\in\bC[z_1,\ldots,z_n]$ be such that
$$F(z_1,\ldots,z_n)\neq 0$$ 
if $\Im(z_k)\ge c_k$, $1\le k\le n$. Then for any $1\le i<j\le n$ one has
$$F\left(z_1,\ldots,z_{i-1},z_i-\frac{\partial}{\partial z_j},z_{i+1},
\ldots,z_j,\ldots,z_n\right)\neq 0$$
whenever $\Im(z_k)\ge c_k$, $1\le k\le n$.
\end{proposition}

From Proposition~\ref{g-LS} we immediately get the following.

\begin{corollary}[(\cite{LS})]\label{LS}
Suppose that $1\leq i < j \leq n$ and 
$F(z_1, \ldots, z_n) \in \HH_n(\CC)$. Then 
$F(z_1,\ldots, z_{i-1}, -\partial_j, z_{i+1}, \ldots,z_j, \ldots, z_n) \in \HH_{n-1}(\CC)\cup\{0\}$.
\end{corollary}

We can now prove the following extension of a famous composition
theorem of Schur \cite{schur} and related results of Mal\'o-Szeg\H{o} 
\cite{CC1,RS} to the multivariate case. Further 
consequences of Theorem \ref{multischur} will be given in \S \ref{s54} and 
\S \ref{ss-polya}. 

\begin{theorem}\label{multischur}
Assume that all the zeros of 
$f(z)=\sum_{i=0}^{r}a_iz^i\in \HH_1(\RR)$ are 
nonpositive and that 
$F(z_1,\ldots,z_n)=\sum_{j=0}^sQ_j(z_2,\ldots,z_n)z_1^j \in \HH_n(\CC)$
and set $m=\min(r,s)$. Then 
$$
\sum_{k=0}^m k! a_k Q_k(z_2,\ldots,z_n)z_1^k \in \HH_n(\CC)\cup\{0\}.
$$
\end{theorem}
 
\begin{proof}
Suppose that $f$ has only nonpositive zeros. Then 
$f(-z_0w_0) \in \HH_2(\RR)$ so that 
$$G(w_0, z_0, \ldots, z_n):= f(-z_0w_0)F(z_1, \ldots, z_n) \in \HH_{n+2}(\CC).$$ 
By Corollary~\ref{LS} we have that 
$$H(z_0, z_1,\ldots, z_n):=
G\!\left(-\frac{\partial}{\partial z_1}, z_0, \ldots, z_n\right) \in \HH_{n+1}(\CC)\cup 
\{0\}.$$ 
This means that 
$$
H(z_1, 0, z_2,\ldots, z_n)= \sum_{k=0}^m k! a_k Q_k(z_2,\ldots,z_n)z_1^k \in 
\HH_n(\CC)\cup\{0\},
$$
as required.
\end{proof}

\subsection{Proofs of Theorems \ref{multichar} and \ref{finitemulti}}

We can now settle the classification of multivariate  multiplier sequences stated in Theorem \ref{multichar}. 

\begin{proof}[of Theorem \ref{multichar}]
For the ``only if'' direction what remains to be proven is the statement 
about the signs. If it were false for some (multivariate) multiplier sequence $\lambda$ then since $\lambda$ is a product of univariate multiplier sequences whose entries either all have the same sign or alternate in sign  there would exist  
$\alpha \in \NN^n$ such that $\lambda(\alpha) \neq 0$ and 
$\lambda(\alpha +e_i)\lambda(\alpha +e_j)<0$.    
Let $T$ be the corresponding (diagonal) operator and apply it 
to $z^\alpha(1-z_iz_j) \in \HH_n(\RR)$. By Lemma \ref{det0} we get 
$$
T(z^\alpha(1-z_iz_j)) = \lambda(\alpha)z^\alpha
\left(1-\frac{\lambda(\alpha +e_i)\lambda(\alpha+e_j)}{\lambda(\alpha)^2}z_iz_j\right).  
$$ 
Since $1 + a z_1z_2 \in \HH_2(\RR)$ if and only if $a \leq 0$ this is a contradiction.

Now $\alpha \mapsto \lambda(\alpha)$ is a multiplier sequence 
if and only if $\alpha \mapsto (-1)^{|\alpha|}\lambda(\alpha)$ is a 
multiplier sequence so we may assume that $\lambda(\alpha) \geq 0$ for 
all $\alpha \in \NN^n$. By applying the $\lambda_i$'s one at a time 
we may further assume that $\lambda_i \equiv 1$ for $2 \leq i \leq n$. Hence 
we have to show that if 
$f(z_1,\ldots,z_n):= \sum_{i= 0}^MQ_i(z_2,\ldots, z_n)z_1^i \in \HH_n(\RR)$ 
and 
$\lambda : \NN \rightarrow \RR$ is a nonnegative univariate  
multiplier sequence then 
$$
\sum_{i=0}^M\lambda(i)Q_i(z_2,\ldots, z_n)z_1^i \in \HH_n(\RR)\cup\{0\}. 
$$
By Theorem \ref{ps} there are polynomials  
$$
p_k(z) = \sum_{i=0}^{N_k}\frac{\lambda_{i,k}}{i!}z^i,\quad k\in\bN,
$$
with only  nonpositive zeros and 
$$
\lim_{k \rightarrow \infty}p_k(z) = \sum_{i=0}^\infty\frac{\lambda(i)}{i!}
z^i. 
$$
Furthermore, by 
Theorem \ref{multischur} we know that 
$$
\sum_{i=0}^M\lambda_{i,k}Q_i(z_2, \ldots, z_n)z_1^i \in 
\HH_n(\RR)\cup\{0\}
$$
for each $k$. Since  $\lim_{k \rightarrow \infty}\lambda_{i,k} = \lambda(i)$ 
we get 
$
\sum_{i=0}^M\lambda(i)Q_i(z_2,\ldots, z_n)z_1^i \in \HH_n(\RR)\cup\{0\}
$
by Lemma \ref{althurwitz}, which settles the theorem.
\end{proof}

Let us finally prove  the characterization of finite order 
multiplier sequences. 

\begin{proof}[of Theorem \ref{finitemulti}]
The ``if'' direction is an immediate consequence of Theorem~\ref{mainreal},  
the proof of which is given in \S \ref{s-main} below. To prove the 
converse statement note first that Theorem \ref{multichar} implies 
that the symbol of $T$ is 
the product of the symbols of the corresponding  univariate operators. 
Hence it suffices to settle the case $n=1$. Let 
$F(z,w) = \sum_{k=0}^Na_kz^kw^k$ be the symbol of $T$. By Theorem 
\ref{ps} we know that all zeros of 
$$
g_m(z)=T[(1+z)^m]= \sum_{k=0}^Ma_k(m)_kz^k(1+z)^{m-k} 
$$
are real and have the same sign. Note that these zeros are actually 
nonpositive since $z=-1$ is a zero of $g_m(z)$ for all 
large $m$. Now 
$$
g_m(z/m)= \left(1+ \frac z m\right)^m \sum_{k=0}^Ma_k\frac{(m)_k}{m^k}z^k(1+z/m)^{-k}
$$
and since $\lim_{m \rightarrow \infty}\frac{(m)_k}{m^k}=1$ we have 
$$
\lim_{m \rightarrow \infty} g_m(z/m) = e^z\sum_{k=0}^Ma_kz^k.
$$
Hence by Lemma \ref{althurwitz} the polynomial $\sum_{k=0}^Ma_kz^k$ has all 
nonpositive zeros and the theorem 
follows.
\end{proof}

\section{Algebraic and geometric properties of stability 
preservers}\label{s-main}

\subsection{Sufficiency in Theorems~\ref{maincomplex} and \ref{mainreal}}

Since $\HH_n(\RR) \subset \HH_n(\CC)$ it is enough to prove only the sufficiency in 
Theorem~\ref{maincomplex}.  Recall the ``affine differential contraction''
of a polynomial  
$F\in \CC[z_1,\ldots,z_n]$ defined in \eqref{not-contr} and note that 
the following consequence of Corollary~\ref{LS} actually settles the 
sufficiency part in Theorem~\ref{maincomplex}.

\begin{corollary}\label{sufficiency}
Let $T \in \A_n[\CC]$ and suppose that 
$F_T(z,-w) \in \HH_{2n}(\CC)$. Then $T \in \A_n(\CC)$. 
\end{corollary}

\begin{proof}
If $F_T(z,-w) \in \HH_{2n}(\CC)$ and 
$f(v) \in \HH_{n}(\CC)$ then $F_T(z,-w)f(v) \in \HH_{3n}(\CC)$. 
By Corollary \ref{LS} if we exchange the variables 
$w_i$'s for $-{\partial}/{\partial v_i}$'s the resulting polynomial 
will be in $\HH_{2n}(\CC)\cup\{0\}$. If we then replace each 
variable $v_i$ with $z_i$ for all  $1\le i\le n$, we get a polynomial in 
$\HH_{n}(\CC)\cup\{0\}$. This polynomial is indeed $T(f)$. 
\end{proof}

\subsection{Necessity in Theorems~\ref{maincomplex} and \ref{mainreal}}\label{necessity}
Let $T = \sum_{\alpha,\beta}a_{\alpha \beta}z^\alpha \partial^\beta \in 
\A_n[\RR]$. We may write $T$ as a finite sum $T = 
\sum_{\gamma \in \ZZ^n}z^\gamma T_\gamma$, where 
$T_\gamma = \sum_{\beta}a_{\gamma+\beta,\beta}z^\beta \partial^\beta$.     
It follows that $T_\gamma$ acts on monomials as 
$T_\gamma(z^\alpha)= \lambda_{\gamma}(\alpha)z^\alpha$ for some 
function $\lambda_\gamma : \NN^n \rightarrow \RR$. The following lemma 
gives a sufficient condition for $\lambda_\gamma$ to be a multiplier 
sequence. 

\begin{lemma}\label{cxhull}
Let $T = \sum_\gamma z^\gamma T_\gamma\in \A_n(\RR)$ and denote by 
$CH(T)$ the convex hull of the set $\{ \gamma : T_\gamma \neq 0\}$. 
If $\kappa \in \ZZ^n$ is a vertex (face of dimension $0$) of $CH(T)$ then 
$T_\kappa$ is a multiplier sequence. 
\end{lemma}
\begin{proof}
Let  $v \in \RR_+^n$. If $f(z) \in \HH_n(\RR)$ then 
$f(vz)= f(v_1z_1, \ldots, v_nz_n) \in \HH_n(\RR)$ and 
$$
T[f(vz)]= \sum_{\alpha, \beta}a_{\alpha \beta}z^\alpha v^\beta 
(\partial^\beta f)(vz).
$$
Hence
$$
T^v := \sum_{\alpha,\beta}a_{\alpha \beta}v^{\alpha-\beta}
z^\alpha \partial^\beta = 
\sum_{\gamma}v^\gamma z^\gamma T_\gamma \in \A_n(\RR). 
$$
Let $\langle z, \mu \rangle = a$ be a supporting hyperplane 
of the vertex $\kappa$. Hence, up to replacing $\mu$ with $-\mu$, if necessary,
we have $\langle \gamma-\kappa, \mu \rangle < 0$ for all 
$\gamma \in CH(T) \setminus \{\kappa\}$. Now let 
$v_i=v_i(t) = e^{\mu_i t}$. Then 
$$
v^{-\kappa}T^v = z^\kappa T_\kappa + 
\sum_{\gamma \neq \kappa}
e^{t\langle \gamma-\kappa, \mu \rangle}z^\gamma T_\gamma.
$$  
By letting $t \rightarrow \infty$ we have that 
$z^\kappa T_\kappa \in \A_n(\RR)$ and the lemma follows. 
\end{proof}

Let $f=\sum_\alpha a(\alpha) z^\alpha \in \RR[z_1,\ldots,z_n]$. Define the {\em support}  
$\supp(f)$ of $f$ to be the set $\{ \alpha \in \NN^n : a(\alpha) \neq 0\}$ 
and let $d = \max\{|\alpha|: \alpha \in \supp(f)\}$. We further define the 
{\em leading part} of $f$ to be $a(\alpha)z^\alpha$, where $\alpha$ is the maximal element of the set 
$\{ \gamma \in \supp(f) : |\gamma|=d \}$ with 
respect to the lexicographic order on $\ZZ^n$. Similarly, if 
$T=\sum_{\gamma}z^\gamma T_\gamma \in \A_n[\RR]$ let 
$k= \max\{ |\alpha|: T_\alpha \neq 0\}$ and let $\kappa_0$ 
be the maximal element of the set $\{ \alpha : |\alpha|=k, T_\alpha \neq 0 \}$
with respect to the lexicographical order. 
Since $\kappa_0$ is a vertex of $CH(T)$ we know that 
$\lambda_{\kappa_0}$ is a multiplier sequence with a finite symbol whenever 
$T \in \A_n(\RR)$. We say that $T_{\kappa_0}$ is the {\em dominating part} 
of $T$.  Note that the dominating part of $fg$ is the product of the dominating parts of 
$f$ and $g$. Moreover, if $\lambda_{\kappa_0}(\alpha) \neq 0$ then the dominating part of $T(f)$ is 
$\lambda_{\kappa_0}(\alpha)a(\alpha)z^{\alpha+\kappa_0}$, where $a(\alpha)z^{\alpha}$ is the dominating part of $f$ and $T_{\kappa_0}$ is 
the dominating part of $T$. 

We are now ready to prove that a real stability preserver also preserves 
proper position. Equivalently, Theorem \ref{tevanster} below asserts that 
$\A_n(\RR) \subset \A_n(\CC)$.

\begin{theorem}\label{tevanster}
Suppose that $T \in \A_n(\RR)$ and that $f, g \in \HH_n(\RR)$ are 
such that $f \ll g$. Then either 
$T(f) \ll T(g)$ or $T(f) =T(g)\equiv 0$.
\end{theorem}

\begin{proof}
Let $T_{\kappa_0}$ be the dominating part of $T$. 
We first assume that $f=\sum_\al a_\al z^\al, g=\sum_\al b_\al z^\al \in \HH_n(\RR)$ are 
such that $f \ll g$,  
$0 \leq \deg(f) < \deg(g)$ and $T_{\kappa_0}(f)T_{\kappa_0}(g) \neq 0$. Let the leading parts
of $f$ and $g$ be $a(\alpha)z^\alpha$ and $b(\beta)z^\beta$, 
respectively. Let 
$$
f_H=\sum_{|\al|=\deg(f)}a_\al z^\al \quad \mbox{ and } \quad g_H=\sum_{|\al|=\deg(g)}b_\al z^\al.
$$
By Hurwitz' theorem $f_H$ and $g_H$ are stable. Moreover all coefficients is $f_H$ (and in $g_H$) have the same sign by \cite[Theorem 6.1]{COSW}. Consider 
$f(vt), g(vt) \in \HH_1(\RR)$ where $v \in \RR_+^n$.  Then 
$\deg g(vt) = \deg f(vt)+1$  
 and the signs of the leading coefficients of 
$g(vt)$ and $f(vt)$ will be the same as the signs of $b(\beta)$ and 
$a(\alpha)$, respectively. Since also $f(vt) \ll g(vt)$ we infer that 
$a(\alpha)b(\beta)>0$.  

Now since $T_{\kappa_0}(f)T_{\kappa_0}(g) \neq 0$ it follows that the leading 
parts of $T(f)$ and $T(g)$ are 
$\lambda_{\kappa_0}(\alpha)a(\alpha)z^{\kappa_0+\alpha}$ and 
$\lambda_{\kappa_0}(\beta)b(\beta)z^{\kappa_0+\beta}$, respectively. 
By Theorem~\ref{nObreschkoff} (the multivariate Obreschkoff theorem)   
we know that either 
$T(f) \ll T(g)$ or $T(g) \ll T(f)$. As pointed out in the paragraph preceding
Theorem~\ref{tevanster} dominating parts are necessarily multivariate
multiplier sequences and so by Theorem \ref{finitemulti} 
we have that $\lambda_{\kappa_0}(\alpha)\lambda_{\kappa_0}(\beta)>0$. From the above discussion it follows that for $v \in \RR_+^n$ we have 
$T(f)(vt) \ll T(g)(vt)$ with $\deg(T(f)(vt)) < \deg(T(g)(vt))$, so that $T(f) \ll T(g)$. 

If $\deg(f)>\deg(g)$ we may simply repeat the arguments using 
$-f$ and $g$.  If 
$\deg(f)=\deg(g)$ we consider $f$ and $g + \epsilon z_1f$ with $\eps>0$. Indeed, 
$\deg(f) < \deg(g + \epsilon z_1f)$ and $f \ll g + \epsilon z_1f$ by 
Lemma \ref{cones}. We then apply the argument of the first case, and obtain the desired conclusion as $\eps \rightarrow 0$. 

Suppose now that  
$T_{\kappa_0}(f)T_{\kappa_0}(g) = 0$. There is nothing to 
prove if $fg\equiv 0$.  Let 
$h_\epsilon(z_1,\ldots,z_n)= (1+\epsilon z_1)^{\xi_1} \cdots (1+\epsilon z_n)^{\xi_n}$ with
$\xi_i\in\bN$, $1\le i\le n$,   
and let $f_\epsilon = h_\epsilon f$ and $g_\epsilon = h_\epsilon g$. 
If $\xi=(\xi_1,\ldots,\xi_n)$ is large enough then 
$T_{\kappa_0}(f_\epsilon)T_{\kappa_0}(g_\epsilon) \neq 0$. The 
theorem follows from Lemma \ref{althurwitz} by letting $\epsilon \rightarrow 0$.   
\end{proof}

For the proof of  necessity  in 
Theorems~\ref{maincomplex} and \ref{mainreal} we need to establish first a key
property for symbols of (real) stability preservers.

\begin{lemma}\label{contraction}
Suppose that $F(z,w) \in \RR[z_1,\ldots, z_n,w_1,\ldots, w_n]$ is the  
symbol of an operator in $\A_n(\RR)$ and let  $\la \in (0,1)^n$. 
Then $F(z,\la w)$ is also the symbol of an operator in $\A_n(\RR)$.  
\end{lemma}

\begin{proof}
Suppose that $T \in \A_n(\RR)$ has symbol 
$F(z_1,\ldots,z_n,w_1,\ldots,w_n)$. 
We claim that if $\delta \geq 0$ then the linear operator  
$\E_1^\delta T$ defined by 
$$
\E_1^\delta T(f)=\sum_{m=0}^\infty \frac{\delta^m z_1^m 
T(\partial_1^m f)}{m!}
$$
is an operator $\E_1^\delta T : \HH_n(\RR) \rightarrow 
\HH_n(\RR)\cup \{0\}$. If 
$T_\delta$ is the linear operator with symbol 
$F(z_1,\ldots,z_n,w_1/(1+\delta),\ldots,w_n)$ then a simple calculation 
shows that 
$$
T_\delta(f)= \E_1^\delta T(f(z_1(1+\delta),\ldots,z_n)).
$$
Hence the claim would prove the lemma. 

In order to prove the remaining claim let $\delta \geq 0$ and 
define a linear operator 
$R_\delta T  : \RR[z_1,\ldots,z_n] \rightarrow \RR[z_1,\ldots,z_n]$ by 
$$
R_\delta T(f) = T(f) + \delta z_1T(\partial_1 f).
$$
Suppose that $f \in \HH_n(\RR)$ and that $T(\partial_1 f) \neq 0$. 
Since $1-iw_1 \in \HH_n(\CC)$ we know by Corollary \ref{sufficiency} that 
$1+i\partial_1 \in \A_n(\CC)$, so $f+i\partial_1f \in \HH_n(\CC)$, i.e.,  $\partial_1 f \ll f$. By 
Theorem \ref{tevanster} we know that  
$T(\partial_1 f) \ll T(f)$ and 
$T(\partial_1 f) \ll 
z_1T(\partial_1 f)$, which  
by Lemma \ref{cones} gives $R_\delta T(f) \in \HH_n(\RR)\cup\{0\}$. 
If $T(\partial_1 f)=0$ then 
$R_\delta T(f)=T(f) \in \HH_n(\RR)\cup\{0\}$, so 
$R_\delta T \in \A_n(\RR)$.

An elementary computation shows that 
when we apply $R_\delta$ to $T$ $m$ times we get 
$$
R^m_\delta T(f) = \sum_{k=0}^m \binom m k 
\delta^kz_1^kT(\partial_1^k f). 
$$
By induction, $R^m_{\delta/m} : \HH_n(\RR) \rightarrow \HH_n(\RR)\cup \{0\}$ for 
all $m\in\bN$. Now  
\begin{eqnarray*}
R^m_{\delta/m} T(f) &=& \sum_{k=0}^m \binom m k m^{-k} 
\delta^kz_1^kT(\partial_1^k f) \\
&=& \sum_{k=0}^m \left(1-\frac 1 m\right)\left(1-\frac 2 m\right)\cdots \left(1-\frac{(k-1)}m\right)
\frac{\delta^k z_1^k T(\partial_1^k f)}{k!}. 
\end{eqnarray*}
It follows that $R^m_{\delta/m} T(f)$ tends uniformly to  $\E^\delta T(f)$ 
on any compact subset of $\CC^n$ as $m\to\infty$. 
Thus $\E^\delta T :  \HH_n(\RR) \rightarrow \HH_n(\RR)\cup \{0\}$ by 
Lemma 
\ref{althurwitz}. 
\end{proof}

From Lemma~\ref{contraction} one can easily see that symbols of (real) 
stability preservers actually satisfy the following homotopical property:

\begin{theorem}\label{t-hom}
If $F(z,w) \in \RR[z_1,\ldots, z_n,w_1,\ldots, w_n]$ is the  
symbol of an operator in $\A_n(\RR)$ then $F(\mu z,\la w)$ is also the 
symbol of an operator in $\A_n(\RR)$ for any 
$(\mu,\la)\in [0,1]^n\times [0,1]^n$. Moreover, the corresponding statement 
holds for symbols of operators in $\A_n(\CC)$.
\end{theorem}

We now have all the tools to accomplish the proof of the necessity  
in Theorem \ref{mainreal}.

\begin{proof}[of Theorem \ref{mainreal}]
The final step in the proof is to show  that $F(z,\mu z^{-1}) \neq 0$ whenever $F$ is the symbol of an operator  
$T \in \A_n(\RR)$, $\mu \in \RR_+^n$ and 
$z =(z_1,\ldots,z_n)\in \CC^n$ is such that $\Im(z_i) >0$ for all $1 \leq i \leq n$. Indeed, since $F(z,w)$ is the symbol of a real stability  
preserver if and only if $F(z+\alpha,w)$ is the symbol of a real stability 
preserver 
for all $\alpha \in \RR^n$ the claim implies that 
$F(z+\alpha,\mu z^{-1}) \neq 0$ whenever $F$ is the symbol of an operator  
$T \in \A_n(\RR)$, $\alpha \in \RR^n$, $\mu \in \RR_+^n$ and 
$z \in \CC^n$ is such that $\Im(z_i) >0$ for $1 \leq i \leq n$. But it is 
straightforward to see that any pair $Z,W \in \CC^n$ such that 
$\Im(Z_i)>0$ and $\Im(W_i)<0$ can be written as 
$Z_i=\alpha_i + z_i$ and $W_i=\mu_iz_i^{-1}$, where $\Im(z_i) >0$, $\al_i\in\bR$ and $\mu_i>0$ for all  
$1 \leq i \leq n$.  Thus, the theorem follows from this claim. 

Let $T = \sum_{\alpha,\beta}a_{\alpha \beta} z^{\alpha}\partial^\beta \in 
\A_n(\RR)$ and let $F$ be its symbol. By multiplying with a large 
monomial we may assume that $a_{\alpha \beta}=0$ if $\alpha \ngeq \beta$.  
Let $v \in \RR_+^n$ and denote by 
${}^vT$ the operator with symbol $F(z,vw)$. 
By Lemma \ref{contraction} we have that 
\begin{eqnarray*}
{}^vT(z^\gamma)z^{-\gamma} &=& 
\sum_{\alpha,\beta}a_{\alpha \beta}v^{\beta}z^{\alpha-\beta}(\gamma)_\beta\\ 
&= &\sum_{\alpha,\beta}a_{\alpha \beta}(v\gamma)^{\beta}z^{\alpha-\beta}
(\gamma)_\beta \gamma^{-\beta}
\in \HH_n(\RR)\cup\{0\}
\end{eqnarray*}
for all $v \in (0,1)^n$. 
Fix $\mu \in \RR_+^n$ and let $v$ in the above equation be of the form 
$\mu \gamma^{-1}$ with $\ga=(\ga_1,\ldots,\ga_n)\in\bN^n$, where
$\ga^{-1}=(\ga_1^{-1},\ldots,\ga_{n}^{-1})$. Then  
$v \in (0,1)^n$ for large $\gamma$. Letting 
$\gamma$ tend to infinity and observing that 
$(\gamma)_\beta \gamma^{-\beta} \rightarrow 1$ we find by Lemma \ref{althurwitz} that  
$$
\sum_{\alpha,\beta}a_{\alpha \beta}\mu^{\beta}z^{\alpha-\beta}
= F(z,\mu z^{-1}) \in \HH_n(\RR)\cup \{0\}. 
$$
We have to prove that $F(z,\mu z^{-1})$ is not identically zero. To do this 
observe that 
$$
F(z,\mu z^{-1})= \sum_{\kappa} z^\kappa 
\sum_{\beta}a_{\beta+\kappa,\beta}\mu^\beta.
$$
By Lemma \ref{cxhull} the dominating part,  
$T_{\kappa_0} = \sum_{\beta}a_{\beta+\kappa_0,\beta}z^\beta \partial^\beta$, of 
$T$ is  an  operator associated to a multiplier sequence with finite symbol.
Hence the nonzero coefficients $a_{\beta+\kappa_0,\beta}$ are all of the 
same sign by Theorem \ref{finitemulti}. This means that the coefficient 
of $z^{\kappa_0}$ in $F(z,\mu z^{-1})$ is nonzero and proves the theorem. 
\end{proof} 

The proof of  necessity in Theorem \ref{maincomplex} now follows easily.
\begin{proof}[of Theorem \ref{maincomplex}]
Let  $T \in \A_n(\CC)$ and write the symbol of $T$ as  $F(z,w) = F_R(z,w)+iF_I(z,w)$, where 
$F_R(z,w)$ and $F_I(z,w)$ have real coefficients. Let further $T_R$ and $T_I$ be the corresponding operators. Now $T: \HH_n(\RR) \rightarrow  \HH_n(\CC)\cup \{ 0 \}$ so by Lemma 
\ref{addz} we have that $T_R+z_{n+1}T_I : \HH_n(\RR) \rightarrow  \HH_{n+1}(\RR)\cup \{ 0 \}$. 
Hence by Lemma \ref{lines} we know that $T_R+(\lambda z_1+ \alpha)T_I \in \A_n(\RR)\cup \{0\}$ for every $\lambda \in \RR_+$ and $\alpha \in \RR$. Suppose that $T_R+(\lambda z_1+ \alpha)T_I =0$. Then $T= (i-\alpha - \lambda z_1)T_I$, so $T_I = 0$ since 
$i-\alpha - \lambda z_1 \notin \HH_1(\CC)$. We thus have $T_R+(\lambda z_1+ \alpha)T_I \in \A_n(\RR)$ for every $\lambda \in \RR_+$ and $\alpha \in \RR$, which by Theorem \ref{mainreal} gives 
$F_R+(\lambda z_1+ \alpha)F_I \in \HH_n(\RR)$ for every $\lambda \in \RR_+$ and $\alpha \in \RR$. By  
Lemma \ref{lines} and Lemma \ref{addz} this implies that $F=F_R + iF_I \in \HH_n(\CC)$, as was 
to be shown. 
\end{proof}

\subsection{The Weyl product and Schur-Mal\'o-Szeg\H{o} type 
theorems}\label{s54}

The results of \S \ref{s32} provide a unifying framework for most
of the classical composition theorems for univariate hyperbolic polynomials
\cite{CC1,M,RS,schur}.
Moreover, they lead to natural multivariate extensions of these
composition theorems. Let us for instance 
consider two operators $S,T\in\A_n[\bC]$ with symbols $F_S(z,w)$ and 
$F_T(z,w)$, respectively. The well-known product formula in the Weyl 
algebra \cite{Bj} asserts that the symbol of the composite operator $ST$ 
is given by
 
\begin{equation}\label{prod-weyl}
F_{ST}(z,w)= \sum_{\kappa \in \NN^n}\frac 1 {\kappa!} \partial^\kappa_wF_S(z,w)
\partial^\kappa_zF_T(z,w).
\end{equation}
This suggests the following definition.

\begin{definition}\label{d-wel}
Let $z=(z_1,\ldots,z_n)$, $w=(w_1,\ldots,w_n)$.
The Weyl product of two polynomials $f(z,w),g(z,w)\in\CC[z,w]$ is given by
$$(f\star g)(z,w)=\sum_{\kappa \in \NN^n}\frac {(-1)^{|\kappa|} } {\kappa !} \partial^\kappa_zf(z,w)
\partial^\kappa_wg(z,w).$$
\end{definition}

Theorems~\ref{maincomplex} and \ref{mainreal} and~\eqref{prod-weyl} 
imply that the Weyl product of polynomials defined above
preserves (real) stability:

\begin{theorem}\label{t-comp}
If $f(z,w)$ and $g(z,w)$ are (real) stable polynomials in the variables 
$z_1, \ldots, z_n$ and 
$w_1, \ldots, w_n$ then their Weyl product $(f\star g)(z,w)$ is also
(real) stable.
\end{theorem}

\begin{example}[(Schur-Mal\'o-Szeg\H{o} theorem)]
Suppose that $S,T\in\A_1[\bR]$ are such that
$F_S(z,w) = f(\lambda^{-1}zw)$ and $F_T(z,w)= g(\lambda z)$, 
where $f\in \HH_1(\bR)$ has only nonpositive zeros, $g\in \HH_1(\bR)$ and 
$\lambda >0$. Then $S,T\in\A_1(\bR)$ by Theorem~\ref{mainreal}. Hence 
$ST\in\A_1(\bR)$ and therefore
$$F_{ST}(z,-w)= \sum_{k\ge 0}\frac 1 {k!} 
z^k f^{(k)}(-\lambda^{-1}zw)g^{(k)}(\lambda z)\in\HH_2(\bR)$$
by~\eqref{prod-weyl}. Letting $w=0$ and $\lambda \rightarrow 0$ it follows 
that 
$$\sum_{k\ge 0} k!  \frac {f^{(k)}(0)}{k!} \frac {g^{(k)}(0)}{k!} z^k
\in \HH_1(\bR),$$
which is a well-known result of Schur \cite{schur}, Mal\'o and Szeg\H{o} 
\cite{CC1,RS}. 
\end{example}

From Theorem \ref{t-comp} one can also recover de Bruijn's 
composition results \cite{db1,db2}. As for composition (or Hadamard-Schur   
convolution)   
theorems in the multivariate case, we should point out that in \cite{hink} 
Hinkkanen obtained such a result for multi-affine polynomials -- i.e., 
multivariate polynomials of degree at most one in each variable -- that are 
nonvanishing when all variables lie in the open unit disk $\DD$. Note that in 
the case of the open upper half-plane Theorem \ref{t-comp} generalizes 
Hinkkanen's composition theorem to arbitrary (not necessarily multi-affine) 
stable polynomials. In fact, by an appropriate conformal transformation one 
can obtain an analog of Theorem \ref{t-comp} for multivariate polynomials of 
arbitrary degrees that are nonvanishing when all variables are in $\DD$, 
thus extending Hinkkanen's convolution theorem.

Finally, by 
using Theorems~\ref{maincomplex} and \ref{mainreal} 
we can derive yet another property
of (real) stability preservers: 

\begin{proposition}\label{prop-coeff}
Let $T\in\A_n[\bC]$ with $F_T(z,w)=\sum_{\al\in\bN^n}Q_{\al}(z)w^{\al}$, 
where as before $z=(z_1,\ldots,z_n)$ and $w=(w_1,\ldots,w_n)$. 
If $T\in\A_n(\bC)$ 
(respectively, $\A_n(\bR)$) then $Q_{\al}(z)\in\HH_n(\bC)\cup\{0\}$ 
(respectively,
$\HH_n(\bR)\cup\{0\}$) for all $\al\in\bN^n$.
\end{proposition}

\begin{proof}
If $T\in\A_n(\bC)$, then 
$F_T(z_1,\ldots,z_n,-w_1,\ldots,-w_n)\in \HH_{2n}(\bC)$ by 
Theorem~\ref{maincomplex}. It follows that for any polynomial 
$P(v_1,\ldots,v_n)\in\HH_n(\bC)$ one has
$$P(v_1,\ldots,v_n)F_T(z_1,\ldots,z_n,-w_1,\ldots,-w_n)\in \HH_{3n}(\bC)$$ 
hence
$$P\!\left(-\frac{\partial}{\partial w_1},\ldots,-\frac{\partial}{\partial w_n}\right)\!F_T(z_1,\ldots,z_n,-w_1,\ldots,-w_n)\in \HH_{2n}(\bC)\cup\{0\}$$
by Corollary~\ref{LS}. Now 
the polynomial 
$P_{\al}(v_1,\ldots,v_n):=v_1^{\al_1}\cdots v_n^{\al_n}$ clearly belongs to
$\HH_n(\bC)$ for any $\al=(\al_1,\ldots,\al_n)\in\bN^n$, so that
by the above one has
\begin{multline*}\al!Q_{\al}(z)=\\
P_{\al}\!\left(-\frac{\partial}{\partial w_1},\ldots,-\frac{\partial}{\partial w_n}\right)\!F_T(z_1,\ldots,z_n,-w_1,\ldots,-w_n)\bigg|_{w_1=\cdots=w_n=0}\!\!\in \HH_{n}(\bC)\cup\{0\},
\end{multline*}
as required. The case when $T\in\A_n(\bR)$ is treated similarly.
\end{proof}

\subsection{Duality, P\'olya's curve theorem and 
generalizations}\label{ss-polya}

Let us first establish the duality property stated in Theorem~\ref{adjoint}. 

\begin{proof}[of Theorem~\ref{adjoint}]
By Theorem \ref{maincomplex} we have that $T \in \A_n(\CC)$ if and only if 
$$G(z,w):=F_T(z,-w) \in \HH_n(\CC).$$ 
But 
$F_{T^*}(z,-w)= \overline{G(-\overline{w}, -\overline{z})} \in \HH_n(\CC)$ 
so the desired conclusion follows from Theorem \ref{maincomplex}. The same 
arguments combined with Theorem~\ref{mainreal} prove the analogous statement 
for $\A_n(\RR)$.
\end{proof}

In view of Theorem~\ref{adjoint} one can both recover known results and 
deduce new ones by a simple dualization procedure, as illustrated in the 
following examples.

\begin{example}[(Hermite-Poulain-Jensen theorem)]
Let $p(z) = \sum_{k=0}^{n}a_kz^k\in \RR[z]\setminus\{0\}$, 
$T=p\!\left({d}/{dz}\right)=\sum_{k=0}^{n}a_k{d^k}/{dz^k}\in 
\A_1[\RR]$, and let $T_p$ be the linear operator on $\bR[z]$ defined by
$T_p(f)(z)=p(z)f(z)$. Then $T^*=T_p$ so by Theorem~\ref{adjoint} one has  
$T \in \A_1(\RR)$ if and only if $T_p \in \A_1(\RR)$, which clearly holds 
if and only if $p \in \HH_1(\RR)$.
\end{example}

\begin{example}
The main result of \cite{ABH} (Theorem 1.4 in {\em op.~cit.})~shows 
that any operator in $\A_1(\bR)$
that commutes with the ``inverted plane differentiation'' operator 
$D_{\sharp}=z^2D$, where $D={d}/{dz}$, is of the form $\al D_{\sharp}^k$
for some $k\in\NN$ and $\al\in\RR$. Therefore, by Theorem~\ref{adjoint} we 
conclude that
any operator in $\A_1(\bR)$ that commutes with the operator 
$zD^2=D_{\sharp}^*$ is of the
form $\al(zD^2)^k$ for some $k\in\NN$ and $\al\in\RR$.
\end{example}

As we will now explain, both Theorem~\ref{mainreal} and 
Theorem~\ref{adjoint} admit natural
geometric interpretations that lead to further interesting consequences. 
For simplicity's sake, we will
only focus on the case $n=1$. 

\begin{definition}\label{pro-I}
Let $f(z,w)\in \bR[z,w]$ be a nonzero polynomial in two variables of (total)
degree $d$ and define the real algebraic curve $\Ga_f$ 
(of degree $d$) by
$$\Ga_f=\{(z,w)\in\bR^2\,:\,f(z,w)=0\}.$$
We say that $f$, or equivalently $\Ga_f$, has the {\em intersection property} 
$(\cI_+)$ if $\Ga_f$ has $d$ real intersection points (counted with 
multiplicities) with any line in $\bR^2$ of the form
$$
w=\al z+\be,\,\text{ where }\al>0,\,\be\in\bR. 
$$
Similarly, we say that $f$ (or $\Ga_f$) has the {\em intersection property} 
$(\cI_{-})$ if $\Ga_f$ has $d$ real intersection points (counted with 
multiplicities) with any line in $\bR^2$ of the form
$$
w=\al z+\be,\,\text{ where }\al<0,\,\be\in\bR.
$$
The {\em symbol curve} of an operator $T\in \A_1[\bR]$ with symbol 
$F_T(z,w)\in\bR[z,w]$ of degree $d$ is the real algebraic curve 
(of degree $d$) given by 
$$\Ga_T=\{(z,w)\in\bR^2\,:\,F_T(z,w)=0\}.$$
\end{definition}

From Lemma~\ref{lines} and Definition~\ref{pro-I} we get:

\begin{corollary}\label{cor-i+}
Let $f$ be a nonzero real polynomial in two variables. 
Then $f\in\HH_2(\bR)$ if and only if $\Ga_f$ has
the intersection property $(\cI_+)$.
\end{corollary}

Therefore, in the univariate case 
Theorem~\ref{mainreal} may be restated as follows.

\begin{corollary}\label{cor-r1}
Let $T\in\A_1[\bR]$. Then $T\in\A_1(\bR)$ if and only if its symbol
curve $\Ga_T$ has the intersection property $(\cI_{-})$.
\end{corollary}

As depicted in Figure~\ref{fig-hpo} below, Corollary \ref{cor-r1} 
essentially allows one
to visualize whether an operator $T\in\A_1[\bR]$ preserves hyperbolicity by
checking whether all lines in $\bR^2$ with negative slope has the required number
of intersection points with $\Ga_T$.

\begin{figure}[!htb]
\centerline{\hbox{\epsfysize=5cm\epsfbox{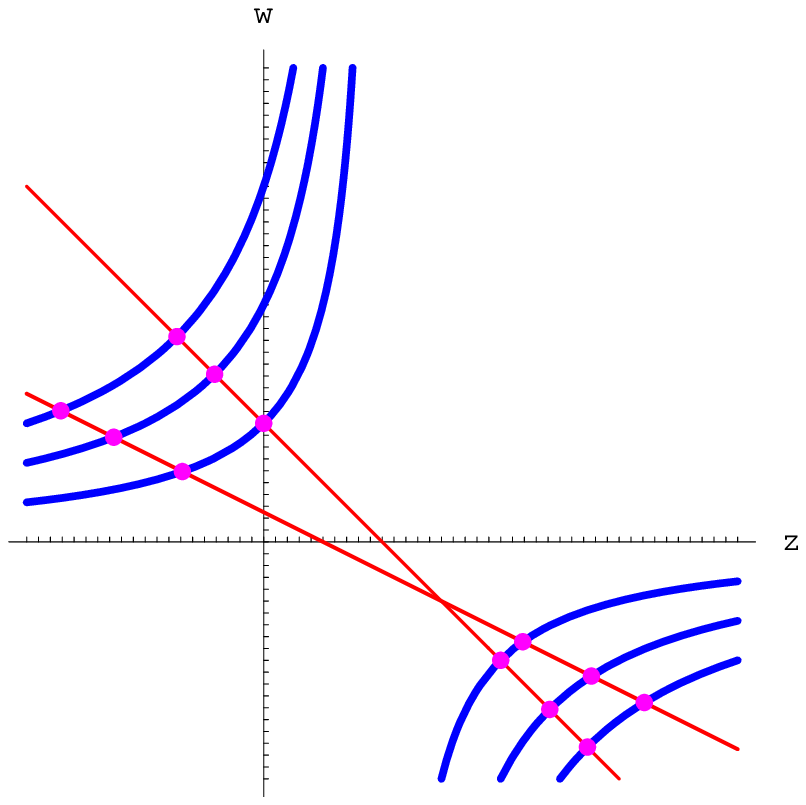}}\hspace{2cm}
\hbox{\epsfysize=5cm\epsfbox{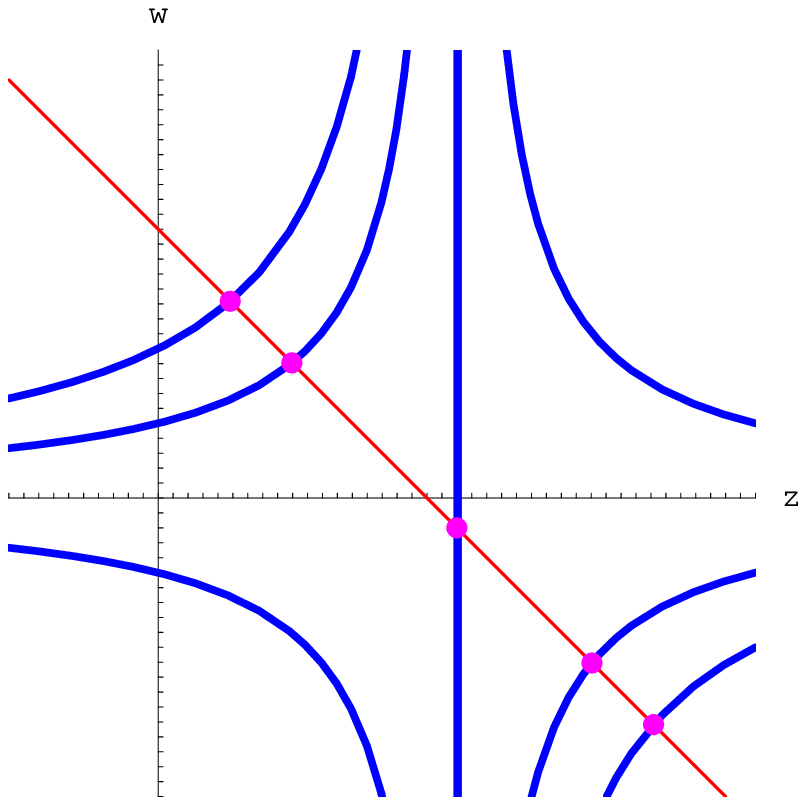}}}
\caption{\label{fig-hpo} Left picture: the symbol curve of degree $d=6$ of 
an operator in 
$\A_1(\bR)$. Right picture: the symbol curve of degree $d=7$ of an operator 
in $\A_1[\bR]$ for which property $(\cI_{-})$ fails.}
\end{figure}

In the same spirit, a simple geometric interpretation and proof of 
Theorem \ref{adjoint} for $n=1$ is as follows: if $T\in\A_1[\bR]$ 
then $F_{T^{*}}(z,w)=F_T(w,z)$ so $\Ga_{T^*}$ is just
the reflection of $\Ga_{T}$ in the main diagonal 
in the $zw$-plane (i.e., the line $w=z$). Since the intersection property 
$(\cI_{-})$ is clearly
invariant under this reflection we conclude that $\Ga_T$ and $\Ga_{T^*}$ 
have the aforementioned property simultaneously.

\begin{figure}[!htb]
\centerline{\hbox{\epsfysize=5cm\epsfbox{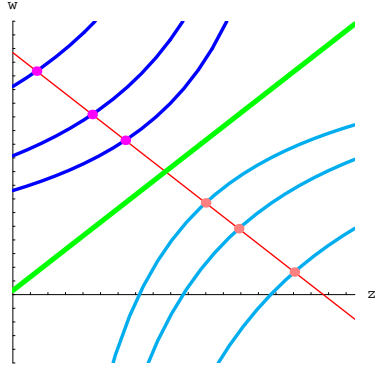}}}
\caption{\label{fig-dual} The symbol curve $\Ga_T$ of degree $d=3$ of an
operator $T\in\A_1(\bR)$ and its dual curve $\Ga_{T^*}$}
\end{figure}

These geometric reformulations of Theorem~\ref{mainreal} and 
Theorem~\ref{adjoint} provide a unifying framework for ``curve type theorems'' 
that include and considerably strengthen P\'olya's original result \cite{P1}
and its various known generalizations \cite{CC3}. 

\begin{example}\label{ex-c1}
In \cite{P1} P\'olya proved the following result that he
considered as the most general 
      theorem on the reality of roots of algebraic equations known at 
      the time (1916).

\begin{theorem}[(P\'olya's curve theorem)]\label{th:polya} 
Let $f(x)$ be a (nonzero) 
	  hyperbolic polynomial of degree $n$, and let 
	  $b_{0}+b_{1}x+\ldots+ b_{n+m}x^{n+m}$, where $m\ge 0$, 
	  be a hyperbolic polynomial with $b_{i}>0$ for all $0\le i\le n$. 
Set
	  $$G_1(x,y)=b_{0}f(y)+b_{1}xf'(y)+b_{2}x^2f''(y)+\ldots+
	  b_{n}f^{(n)}(y)\in\bR[x,y].$$
Then $G$ has the intersection property $(\cI_{+})$.
\end{theorem}

\begin{proof}
By assumption the polynomial $q(x)=\sum_{k=0}^{n+m}b_kx^k$ has only real and
nonpositive zeros hence $q(-xz)\in\HH_2(\bR)$ and thus the 
polynomial in variables $x,y,z$ given by $q(-xz)f(y)$ belongs to $\HH_3(\bR)$. 
From Corollary \ref{LS} we then get $q(xD_y)f(y)\in \HH_2(\bR)$, where
$D_y={d}/{dy}$, and the result follows by Corollary \ref{cor-i+}.
\end{proof} 
\end{example}
	 
\section{Strict stability and strict real stability preservers}\label{s-strict}

A natural question in the present context is to characterize all finite order linear differential operators that 
preserve {\em strict stability}
and {\em strict real stability}, respectively. These notions are defined as 
follows: note first that 
the set of real stable univariate polynomials coincides with
the set of hyperbolic univariate polynomials. Denote
by $\HH_1^{s}(\RR)$ the set of all {\em strictly hyperbolic} univariate 
polynomials, i.e., polynomials in $\HH_1(\RR)$ with only simple zeros.

\begin{definition}\label{d-strict}
A polynomial $f\in \HH_n(\bR)$ is called {\em strictly real stable} if 
$f(\al+vt)\in \HH_1^{s}(\RR)$ for any $\al\in\RR^n$ and $v\in\RR_{+}^n$.
One calls a polynomial $g\in  \HH_n(\bC)$ {\em strictly stable} if 
$g(z_1, \ldots, z_n) \neq 0$ for all $n$-tuples 
$(z_1, \ldots, z_n) \in \CC^n$ with $\Im(z_j) \ge 0$.
\end{definition}

Let $\HH_n^s(\bC)$ (respectively, $\HH_n^s(\bR)$) be the set of all strictly 
stable (respectively, strictly real stable) polynomials in $n$ variables. 
Clearly, if $n\ge 2$ then 
$\HH_n^s(\bR)=\HH_n^s(\bC)\cap \bR[z_1,\ldots,z_n]$ while 
$\HH_1^s(\bC)\cap \bR[z]=\bR\setminus \{0\}$.
Denote by $\A_n^s(\bC)$
and $\A_n^s(\bR)$ the submonoids of $\A_n[\bC]$ and $\A_n[\bR]$ 
consisting of all strict stability and strict real stability
preservers, respectively, i.e., $\A_n^s(\bC)=\{T\in \A_n[\bC]:T(\HH_n^s(\bC))
\subseteq \HH_n^s(\bC)\cup\{0\}\}$ and 
$\A_n^s(\bR)=\{T\in \A_n[\bR]:T(\HH_n^s(\bR))
\subseteq \HH_n^s(\bR)\cup\{0\}\}$. 

In this section we give necessary and sufficient 
conditions in order for a linear
operator to belong to either $\A_n^s(\bC)$
or $\A_n^s(\bR)$.

\begin{theorem}\label{s-maincomplex}
Let $T \in \A_n[\CC]$. If $T \in \A_n^s(\CC)$ then  
$F_T(z_1,\ldots,z_n,-w_1,\ldots,-w_n)\neq 0$ 
whenever $\Im(z_j)\ge 0$ and $\Im(w_k)>0$ for all $1\le j,k\le n$.
\end{theorem}

\begin{theorem}\label{s-mainreal}
Let $T \in \A_n[\RR]$. If $T \in \A_n^s(\RR)$ then  
$F_T(z_1,\ldots,z_n,-w_1,\ldots,-w_n)\neq 0$ 
whenever $\Im(z_j)\ge 0$ and $\Im(w_k)>0$ for all  $1\le j,k\le n$.
\end{theorem}

To prove Theorems~\ref{s-maincomplex} and \ref{s-mainreal} we need to establish
a multivariate extension of the following classical result 
\cite{obreschkoff,RS}, compare with 
Theorems \ref{1Obreschkoff} and \ref{nObreschkoff}.

\begin{theorem}[(Strict Obreschkoff theorem)]\label{strictObreschkoff}
Let $f,g \in \RR[z]$. Then
$$ 
\{\alpha f + \beta g : \alpha, \beta \in \RR, \alpha^2 + \beta^2 >0\} \subset \HH_1^s(\RR)
$$
if and only if $f+ig \in \HH_1^s(\CC)$ or $g+if \in \HH_1^s(\CC)$. 
\end{theorem} 

\begin{theorem}\label{strictnObreschkoff}
Let $f,g \in \RR[z_1,\ldots,z_n]$. Then
$$ 
\{\alpha f + \beta g : \alpha, \beta \in \RR, \alpha^2 + \beta^2 >0\} \subset \HH_n^s(\RR)
$$
if and only if $f+ig \in \HH_n^s(\CC)$ or $g+if \in \HH_n^s(\CC)$. 
\end{theorem} 

\begin{proof}
This is an immediate consequence of Definition~\ref{d-strict} and
Theorem~\ref{strictObreschkoff}.
\end{proof}

\begin{proof}[of Theorem~\ref{s-maincomplex}]
Suppose that $T\in\A_n^s(\bC)$ and let $f \in \HH_n(\CC)$. Then 
$$f_\epsilon(z_1,\ldots,z_n):=f(z_1+i\epsilon, \ldots, z_n+i\epsilon) 
\in \HH_n^s(\CC)$$ 
for all 
$\epsilon >0$, so 
$T(f_\epsilon) \in \HH_n^s(\CC)$. Letting $\eps\to 0$ it follows
from Hurwitz' theorem that 
$T(f) \in \HH_n(\CC) \cup \{0\}$ and thus 
$T\in\A_n(\bC)$. Now 
Theorem~\ref{maincomplex} implies that 
$$F_T(z_1,\ldots,z_n,-w_1,\ldots,-w_n)
:=\sum_{\al\in\bN^n}Q_{\al}(z_1,\ldots,z_n)(-w)^{\al}\in\HH_{2n}(\bC),$$
where $Q_{\al}\in\bC[z_1,\ldots,z_n]$ are not identically zero only for a 
finite number of multi-indices $\al\in\bN^n$. Fix 
$(z_1^0,\ldots,z_n^0)\in\bC^n$ with $\Im(z_i^0)\ge 0$, $1\le i\le n$. 
Then for any
$\eps>0$ one has $\Im(z_k^0+i\eps)>0$ for all $1\le k\le n$. Hence 
$$\sum_{\al\in\bN^n}Q_{\al}(z_1^0+i\eps,\ldots,z_n^0+i\eps)(-w)^{\al}
\in\HH_{n}(\bC)$$
and by letting $\eps\to 0$ we deduce that
$$\sum_{\al\in\bN^n}Q_{\al}(z_1^0,\ldots,z_n^0)(-w)^{\al}\in\HH_{n}(\bC)\cup
\{0\}.$$
If $\sum_{\al\in\bN^n}Q_{\al}(z_1^0,\ldots,z_n^0)(-w)^{\al}\equiv 0$ then
$Q_{\al}(z_1^0,\ldots,z_n^0)=0$ for all $\al\in\bN^n$ and consequently 
$T(f)(z_1^0, \ldots, z_n^0) =0$ for 
all polynomials $f$, which contradicts the assumption that 
$T\in\A_n^s(\bC)$ and $\Im(z_i^0) \geq 0$ for all $1 \leq i \leq n$.
\end{proof}

\begin{proof}[of Theorem~\ref{s-mainreal}]
Suppose that $T \in \A_n^s(\RR)$. If the symbol is not as in the statement of the theorem then by arguing as in the proof of the necessity in Theorem~\ref{s-maincomplex} we see that there 
exists $Z^0=(z_1^0, \ldots, z_n^0)\in \CC^n$ with $\Im(z_i^0) \geq 0$, $1\le i\le n$, and 
\begin{equation}\label{noll}
F(Z^0, -w) = 0
\end{equation}
as a polynomial in $w=(w_1, \ldots, w_n)$.  Choose a polynomial $f \in \HH_n^s(\RR)$ of sufficiently high degree 
so that $T(\alpha f +\beta \partial_1f) \neq 0$ whenever $\alpha^2+\beta^2 >0$. By 
Theorem~\ref{strictnObreschkoff} we have $T(f)+i T(\partial_1f) \in \HH_n^s(\CC)$. This is however a contradiction since by \eqref{noll} we have 
$
T(f)(Z^0)+iT(\partial_1f)(Z^0) = 0. 
$
\end{proof}

The next two theorems give sufficient conditions for operators in the Weyl 
algebra to be 
strict stability or strict real stability preserving, respectively.

\begin{theorem}\label{suff-ss}
Let $T\in\A_n[\bC]$. If $F_T(z_1,\ldots,z_n,-w_1,\ldots,-w_n)\in\HH_{2n}^{s}(\bC)$ then $T\in\A_n^{s}(\bC)$.
\end{theorem}

\begin{theorem}\label{suff-srs}
Let $T\in\A_n[\bR]$. If $F_T(z_1,\ldots,z_n,-w_1,\ldots,-w_n)\in\HH_{2n}^{s}(\bR)$ then $T\in\A_n^{s}(\bR)$.
\end{theorem}

\begin{proof}[of Theorem~\ref{suff-ss}]
Let $T\in\A_n[\bC]$ and suppose that
$$F_T(z_1,\ldots,z_n,-w_1,\ldots,-w_n)\neq 0$$
whenever $\Im(z_i)\ge 0$ and $\Im(w_j)\ge0$ for all  $1\le i,j\le n$. If
$f(v_1,\ldots,v_n)\in\HH_n^s(\bC)$ then
$$F_T(z_1,\ldots,z_n,-w_1,\ldots,-w_n)f(v_1,\ldots,v_n)\neq 0$$
provided that $\Im(z_i)\ge 0$, $\Im(w_j)\ge 0$ and $\Im(v_k)\geq 0$ for  all 
$1\le i,j,k\le n$. By Proposition \ref{g-LS} we may replace each variable 
$w_j$ with $w_j-{\partial}/{\partial v_j}$ for all  $1\le j\le n$, to get
$$F_T\!\left(z_1,\ldots,z_n,\frac{\partial}{\partial v_1}-w_1,\ldots,
\frac{\partial}{\partial v_n}-w_n\right)f(v_1,\ldots,v_n)\neq 0$$
whenever $\Im(z_i)\ge 0$, $\Im(v_i)\ge 0$, $\Im(w_i)\ge 0$ for $1\le i \le n$. If we now exchange
each variable $v_j$ by $z_j$ for all $1\le j\le n$, and let $w_i=0$ for all $1\le i\le n$,
 we see that
$$
T(f)(z_1,\ldots,z_n)=F_T\!\left(z_1,\ldots,z_n,
\frac{\partial}{\partial z_1},\ldots,
\frac{\partial}{\partial z_n}\right)f(z_1,\ldots,z_n)\neq 0$$
whenever $\Im(z_i)\ge 0$, $1\le i\le n$, hence $T(f)\in \HH_n^s(\bC)$. 
\end{proof}

\begin{proof}[of Theorem~\ref{suff-srs}]
If $F_T$ is as in the statement of the theorem then by Theorem~\ref{suff-ss} we have that 
$T \in \A_n^s(\CC)$. Consider $f \in \HH_n^s(\RR)$. The case when $f$ is a 
nonzero constant, say $f(z)\equiv c\in\bR\setminus \{0\}$, is immediate
since $T(f)(z_1,\ldots,z_n)=cF_T(z_1,\ldots,z_n,0,\ldots,0)\neq 0$ whenever
$\Im(z_i)\ge 0$, $1\le i\le n$, hence $T(f)\in \HH_n^s(\RR)$. Suppose that $f$ is not a 
constant polynomial. By re-indexing the variables  
we may assume that $\partial_1f \not\equiv 0$. Then 
$f+i\partial_1f \in \HH_n^s(\CC)$, so $T(f)+iT(\partial_1 f) \in \HH_n^s(\CC)\cup \{0\}$. By 
Theorem~\ref{strictnObreschkoff} we have that $T(f) \in \HH_n^s(\RR)\cup\{0\}$,
as required.
\end{proof}

To close this section we note that in general 
the necessary conditions stated in 
Theorems~\ref{s-maincomplex} and \ref{s-mainreal} are not sufficient while the
sufficient conditions given in Theorems~\ref{suff-ss}--\ref{suff-srs} are not 
necessary. This may be seen already in the univariate case from the 
following simple examples. The operator $T={d}/{dz}$ is clearly strict
(real) stability preserving but its symbol $F_T(z,w)=w$ does not satisfy the
sufficient conditions stated in Theorems~\ref{suff-ss} and \ref{suff-srs}. 
Consider now the operator $S=2z+1+(z^2+z){d}/{dz}$. One can easily check 
that $F_S(z,-w)\neq 0$ whenever $\Im(z)\ge 0$ and $\Im(w)>0$ and also that $S$ 
preserves strictly real stable (i.e., strictly hyperbolic) polynomials.
However, $S(1)=2z+1\notin\HH_1^{s}(\bC)$ so $S\notin\A_1^{s}(\bC)$. A 
characterization of strict (real) stability preservers would therefore require 
conditions that are intermediate between those of 
Theorems~\ref{s-maincomplex} and \ref{s-mainreal} and 
Theorems~\ref{suff-ss} and \ref{suff-srs}.

\section{Multivariate matrix pencils and applications}\label{s-ap}

We will now give several examples and applications of the above 
results. First we prove Proposition~\ref{pencil} claiming that the polynomial 
\begin{equation*}
f(z_1,\dots, z_n):=\det\left( \sum_{i=1}^nz_iA_i+ B\right)
\end{equation*}
 with  $A_1, \ldots, A_n$  positive semidefinite  matrices and  $B$ a
Hermitian matrix of the same order  is either real stable or 
identically zero.

\begin{proof}[of Proposition~\ref{pencil}]
By a standard continuity argument using Hurwitz' theorem it suffices to prove 
the result only in the case when  all matrices $A_1, A_2,\ldots, A_n$ are 
positive definite. Set $z(t)=\alpha+\la t$ with
$\alpha\in \RR^n$, $\la\in \RR_+^n$ and $t\in \bR$. Note that 
$P:=\la_1A_1+\ldots+\la_n A_n$
is positive definite and thus it has a square root. Then 
$
f(z(t))=\det(P)\det(tI+P^{-1/2}HP^{-1/2}),
$
where $H:=B+\al_1A_1+\ldots +\al_n A_n$. Since 
$f(z(t))$ is a constant multiple of the characteristic polynomial of the 
Hermitian matrix $H$, it has only real zeros.
\end{proof}

\subsection{A stable multivariate extension of the Cauchy-Poincar\'e 
theorem}\label{s51} 
Let $A$ be any $n \times n$ complex matrix and define a polynomial 
$C(A,z)= \det(Z-A) \in \CC[z_1, \ldots, z_n]$, 
where $z=(z_1,\ldots,z_n)$ and 
$Z$ is the (diagonal) matrix with entries $Z_{ij} = z_i\delta_{ij}$. 
Given $1 \leq i,j \leq n$ let $A^{ij}$ be the submatrix of $A$ obtained by
deleting row $i$ and column $j$ and 
set $C_{ij}(A,z) = \det\left( (Z-A)^{ij} \right)$. For 
$z=(z_1,\ldots,z_n)$ and $1 \leq i \leq n$ let 
$z\setminus z_i = (z_1, \dots, z_{i-1}, z_{i+1}, \ldots, z_n)$, so that
$C_{ii}(A,z)= C(A^{ii}, z\setminus z_i )$.

\begin{lemma}\label{l-simple}
For $1\le j\le n$ one has $T_j:=1+i{\partial}/{\partial z_j}\in \A_n(\bC)$.
\end{lemma} 
 
\begin{proof}
The symbol $F_{T_j}(z,w)$ of $T_j$ is $F_{T_j}(z,-w)=1-iw_j$ 
and the latter polynomial is stable since it is obviously nonvanishing if 
$\Im(w_j)>0$. The assertion now follows from Theorem~\ref{maincomplex}.
\end{proof}

\begin{theorem}\label{C-P}
If $A$ is a complex Hermitian $n \times n$ matrix then $C(A,z) \in \HH_n(\RR)$
and 
$$
C(A^{jj}, z\setminus z_j) \ll C(A, z)
$$
for $1 \leq j \leq n$.
\end{theorem}

\begin{proof}
Note that since $A$ is Hermitian $C(A,z)$ is real stable by 
Proposition~\ref{pencil}. Now
$$C(A^{jj},z\setminus z_j)=C_{jj}(A,z)=\frac{\partial}{\partial z_j}C(A,z)
\ll C(A,z)$$
by Lemma~\ref{l-simple} and Corollary~\ref{hhC}.
\end{proof}

The above theorem generalizes the classical Cauchy-Poincar\'e theorem stating
that the eigenvalues of a Hermitian matrix and those of any of 
its degeneracy one
principal submatrices interlace. 

An alternative proof of Theorem~\ref{C-P} may
be obtained by using the following consequence of the Christoffel-Darboux
identity \cite{Go}:

\begin{lemma}
Let $A$ be any $n \times n$ matrix with $n \geq 2$ and let $1 \leq i,j \leq n$. Then 
\begin{equation}\label{kartoffel}
C(A,y)C_{ij}(A,x)-C(A,x)C_{ij}(A,y) = \sum_{k=1}^n(y_k-x_k)C_{ik}(A,x)C_{kj}(A,y) .
\end{equation}
\end{lemma}
\begin{proof}
Let $X = (x_i \delta_{ij})$ and $Y=(x_i \delta_{ij})$.  The identity 
\begin{equation*}
(X-A)^{-1}-(Y-A)^{-1}=(X-A)^{-1}(Y-X)(Y-A)^{-1} 
\end{equation*}
obtains by multiplying on the left with $(X-A)$ and on the right with $(Y-A)$. Taking 
the $ij$-th entry on both sides in the above identity and multiplying by $C(A,x)C(A,y)$ yields formula \eqref{kartoffel}. 
\end{proof}

Now let $A$ be a complex Hermitian $n \times n$ matrix with $n \geq 2$ and let $y=z$, $x=\bar{z}$ and $i=j$ in \eqref{kartoffel}. Note that $C_{ij}(A,\bar{z})=\overline{C_{ji}(A,z)}$ and since 
$C_{ii}(A,z)= C(A^{ii}, z\setminus z_i )$ 
we get
\begin{equation}\label{snyggt}
\Im( C(A,z)C(A^{ii},\overline{z\setminus z_i} )) = \sum_{k=0}^n \Im(z_k)|C_{ik}(A,z)|^2.
\end{equation}

Theorem~\ref{C-P} is obviously true for $n=1$ and the general case follows by induction on 
$n$. Indeed, let $\Im(z_j) >0$ for $1 \leq j \leq n$, where $n \geq 2$. By the induction hypothesis we have 
$C(A^{ii},z\setminus{z_i}) \in \HH_{n-1}(\CC)$ and 
then from \eqref{snyggt} we deduce that 
$$
\Im\left( \frac {C(A,z)}{C(A^{ii},z\setminus{z_i})}\right) = \Im\left(\frac {C(A,z)C(A^{ii},\overline{z\setminus z_i})}{|C(A^{ii},z\setminus{z_i})|^2} \right)\geq 
\Im(z_i)>0.
$$
Hence $C(A,z) \neq 0$ and the desired conclusion follows from Lemma \ref{addz}.
   
\subsection{Lax conjecture for real stable polynomials in two variables}\label{s52} 
Here we will prove that all real stable polynomials in two variables $x, y$ can be written as 
$\pm\det(xA+yB+C)$, where  $A$ and $B$ are positive semidefinite (PSD) matrices and $C$ is a symmetric matrix. The 
proof relies on the Lax conjecture that was recently settled in \cite{LPR} by using in an essential way the results of \cite{HV}. 

\begin{theorem}[(\cite{HV,LPR})]\label{lax}
A homogeneous polynomial $p\in\bR[x,y,z]$ is hyperbolic of degree $d$ with respect to the vector 
$e=(1,0,0)$  if and only if there exist two symmetric $d \times d$ matrices $B,C$  such that 
$$
p(x,y,z)=p(e)\det(xI+yB+zC).
$$
\end{theorem}

We will also need Proposition \ref{pro-g} that we proceed to prove.

\begin{proof}[of Proposition \ref{pro-g}] If $f\in\bR[z_1,\ldots,z_n]$
is of degree $d$ then its {\em homogenization} -- 
i.e., the unique homogeneous 
polynomial $f_{H}\in\bR[z_1,\ldots, z_{n+1}]$ of degree $d$ such that 
$f_{H}(z_1, \ldots, z_n,1)=f(z_1,\ldots,z_n)$ -- is simply 
$$f_{H}(z_1,\ldots, z_{n+1})
=z_{n+1}^{d}f(z_1z_{n+1}^{-1},\ldots,z_nz_{n+1}^{-1}).$$
If $f_H$ is hyperbolic with respect to every vector $v\in\bR^{n+1}$ 
such that
$v_{n+1}=0$ and $v_i>0$, for all $1\le i\le n$, it follows in particular that 
for any $\al=(\al_1,\ldots,\al_n)\in\bR^n$ and $(v_1,\ldots,v_n)\in\bR_{+}^n$
the univariate (real) polynomial in $t$ given by 
$$f_H((\al_1,\ldots,\al_n,1)+(v_1,\ldots,v_n,0)t)=f(\al+vt)$$
is not identically zero (since $\lim_{t\to\infty}t^{-d}f(\al+vt)
=f_H(v_1,\ldots,v_n,0)\neq 0$) and has only real zeros. Hence it
belongs to $\HH_1(\bR)$. Thus $f\in\HH_n(\bR)$ by Lemma \ref{lines}.

Conversely, suppose that $f\in\HH_n(\bR)$ has degree $d$ and is given by
$$f(z)=\sum_{\kappa\in\bN^n}a_{\kappa}z^{\kappa},\quad z=(z_1,\ldots,z_n).$$ 
Let 
$\al=(\al_1,\ldots,\al_{n+1})\in\bR^{n+1}$ and $v=(v_1,\ldots,v_{n+1})\in
\bR^{n+1}$ with $v_{n+1}=0$ and  $v_i>0$ for all  $1\le i\le n$. Since $a_{\kappa}\neq 0$ 
for some $\kappa\in\bN^n$ with $|\kappa|=d$, Hurwitz'
theorem yields
$$g(z):=\lim_{t\to\infty}t^{-d}f(tz)
=\sum_{\kappa\in\bN^n,\,|\kappa|=d}a_{\kappa}z^d\in\HH_n(\bR).$$
Moreover, $g$ is a homogeneous polynomial so by the ``same phase property'' 
established in \cite[Theorem 6.1]{COSW} all nonzero $a_{\kappa}$'s with $|\kappa|=d$ have 
the same sign. Therefore
$$f_H(v)=g(v_1,\ldots,v_n)=\sum_{\kappa\in\bN^n,\,|\kappa|=d}a_{\kappa}
v_1^{\kappa_1}\cdots v_{n}^{\kappa_n}\neq 0$$
since $v_i>0$ for all  $1\le i\le n$. Now, if $\al_{n+1}=0$ then the univariate 
polynomial 
$$t\mapsto f_H(\al+vt)=g(\al_1+v_1t,\ldots,\al_n+tv_n)$$
has only real zeros by Lemma \ref{lines}, while if 
$\al_{n+1}>0$ then again by Lemma \ref{lines} the univariate 
polynomial
$$t\mapsto 
f_H(\al+vt)=\al_{n+1}^{d}f(\al_1\al_{n+1}^{-1}+v_1\al_{n+1}^{-1}t,\ldots,
\al_n\al_{n+1}^{-1}+v_n\al_{n+1}^{-1}t)$$
has only real zeros. By the last part of Lemma \ref{lines},
the same holds when $\al_{n+1}<0$. Hence $f_H$ is hyperbolic with respect to
all vectors $v\in\bR^{n+1}$ as above.
\end{proof}

\begin{lemma}\label{dia}
Suppose that $p\in\bR[x,y,z]$ is a homogeneous polynomial of degree $d$ which is hyperbolic with respect to 
any $(v_1,v_2,0)\in\bR^3$ with $v_1,v_2 \in \RR_+$. Then   
$$
p(x,y,0)= \sum_{i=0}^d a_i x^{d-i}y^i , 
$$
where the $a_i$'s are such that $\sum_{i=0}^d a_i t^i$ is a polynomial with only nonpositive 
zeros. 
\end{lemma}

\begin{proof}
By letting $z \rightarrow 0$ it follows from Hurwitz' theorem 
that $p(x,y,0)$ is hyperbolic with respect to all 
$v \in \RR_+^2$, so it is stable by Proposition \ref{pro-g}. By \cite[Theorem 6.1]{COSW} all the coefficients have the same sign.  Since $p(x,y,0)$ is stable all the zeros of the  polynomial  
where $p(1,t,0)=a_0 + a_1t + \cdots +a_d t^d$ are real and since the coefficients have the same sign all the zeros are nonpositive.
\end{proof}

\begin{theorem}\label{lax-2}
A homogeneous polynomial $p\in\bR[x,y,z]$ of degree $d$ is hyperbolic with respect to all vectors of the form 
$(v_1,v_2,0)$ with $v_1,v_2 \in \RR_+$ if and only if there exist two positive semidefinite 
$d \times d$ matrices $A$ and $B$ and a  symmetric $d \times d$ matrix $C$  such that 
$$
p(x,y,z)=\al\det(xA+yB+zC),
$$
where $\al\in\bR$. Moreover, $A$ and $B$ can be chosen so that $A+B=I$. 
\end{theorem}

\begin{proof}
Let $p$ be hyperbolic of degree $d$ with respect to all vectors of the form 
$(v_1,v_2,0)$ with $v_1,v_2 \in \RR_+$ and let $\al:=p(1,1,0)\in\bR\setminus\{0\}$. 
Consider the polynomial $f(x,y,z)=p(x,x+y,z)$. Then $f(x,y,z)$ is hyperbolic of degree $d$ with respect to all vectors $(v_1,v_2,0)$, where $v_1,v_2 \in \RR_+$. Moreover, it is hyperbolic with respect 
to the vector $e=(1,0,0)$. Hence by Theorem \ref{lax} there exist two symmetric $d \times d$ matrices $B,C$  such that 
$$
f(x,y,z)=f(e)\det(xI+yB+zC).
$$
Since $f$ is hyperbolic with respect to all vectors of the form 
$(v_1,v_2,0)$ with $v_1,v_2 \in \RR_+$ we know by Lemma \ref{dia} that all the eigenvalues 
of $B$ are nonnegative. Hence $B$ is a PSD matrix. Let $A=I-B$. Then 
$$
p(x,y,z)= \al\det(xA+y(I-A)+zC), 
$$
and by Lemma \ref{dia} all zeros of the polynomial 
$$
r(t):=\al^{-1}p(1,t,0)=(1-t)^d\det\left(A+\frac{t}{1-t}I\right)\in\bR[t]
$$
are nonpositive. Inverting this we have 
$$
\det(A+tI)= (1+t)^dr\left(\frac{t}{1+t}\right), 
$$
which implies that $A$ has only nonnegative eigenvalues, so that $A$ is a 
PSD matrix. 
\end{proof}

From Theorem~\ref{lax-2} and Proposition~\ref{pro-g} we deduce the following 
converse to Proposition~\ref{pencil} for real stable polynomials in 
two variables.

\begin{corollary}\label{cor-2var}
Let $f(x,y)\in\bR[x,y]$ be of degree $n$. Then $f$ is real stable if and only 
if there exist two $n\times n$ PSD matrices $A,B$ and a symmetric $n\times n$ 
matrix $C$ such that
$$f(x,y)=\pm\det(xA+yB+C).$$
\end{corollary}

\subsection{Hyperbolicity preservers via determinants and homogenized
symbols}\label{s53} 
Using Theorem \ref{mainreal} with $n=1$ and Corollary~\ref{cor-2var}
we immediately get the following determinantal description of finite order
linear preservers of univariate real stable (i.e., hyperbolic) polynomials.

\begin{theorem}\label{crit1}
Let $T\in\A_1[\bR]$. Then $T\in\A_1(\bR)$ if and only if there exist 
$\al\in\bR$, $d\in\bN$, 
two positive semidefinite 
$d \times d$ matrices $A$ and $B$ and a symmetric $d \times d$ matrix $C$  
such that 
$$T=\al\det(zA-wB+C)\Big|_{w={\partial}/{\partial z}}.$$
\end{theorem}

From Theorem~\ref{crit1} and Proposition~\ref{pro-g} we deduce yet another
characterization of univariate hyperbolicity preservers involving 
real homogeneous (G\aa rding) hyperbolic polynomials in 3 variables:

\begin{theorem}\label{crit2}
Let $T\in\A_1[\bR]$ with symbol $F_T(z,w)$ of 
degree $d$ and let $\tilde{F}_T(y,z,w)$ be the (unique) homogeneous degree $d$
polynomial such that $\tilde{F}_T(1,z,w)=F_T(z,w)$. Then 
$T\in \A_1(\bR)$ if and only if the following conditions hold:
\begin{enumerate}
\item[(i)] $\tilde{F}_T(y,z,w)$ is hyperbolic with respect to $(0,1,1)$,
\item[(ii)] all zeros of $\tilde{F}_T(0,t,1)$ lie in $(-\infty,0]$.
\end{enumerate}
\end{theorem}

\end{document}